\documentstyle[amssymb,amsfonts]{amsart}

\def\Energy{{{\operatorname{E}}}}
\def\R{{{\mathbf R}}}
\def\C{{{\mathbf C}}}

\def\Z{{{\mathbf Z}}}

\def\K{{\mathbf{K}}}
\def\B{{\mathbf{B}}}

\def\N{{\mathcal N}}
\def\eps{\varepsilon}
\newenvironment{proof}{\noindent {\bf Proof} }{\endprf\par}
\def \endprf{\hfill  {\vrule height6pt width6pt depth0pt}\medskip}
\def\emph#1{{\it #1}}
\def\textbf#1{{\bf #1}}

\theoremstyle{plain}
  \newtheorem{theorem}[subsection]{Theorem}

  \newtheorem{proposition}[subsection]{Proposition}
  
  \newtheorem{lemma}[subsection]{Lemma}
  \newtheorem{corollary}[subsection]{Corollary}

\theoremstyle{remark}
  \newtheorem{remark}[subsection]{Remark}
  \newtheorem{remarks}[subsection]{Remarks}
  \newtheorem{example}[subsection]{Example}
  \newtheorem{examples}[subsection]{Examples}

\theoremstyle{definition}
  \newtheorem{definition}[subsection]{Definition}

\include{psfig}

\begin{document}

\title{Product set estimates for non-commutative groups}

\author{Terence Tao}
\address{Department of Mathematics, UCLA, Los Angeles CA 90095-1555}
\email{tao@@math.ucla.edu}
\thanks{T. Tao is supported by a grant from the Packard Foundation.}

\begin{abstract}
We develop the Pl\"unnecke-Ruzsa and Balog-Szemer\'edi-Gowers theory of sum set estimates in the non-commutative setting, with discrete, continuous, and metric entropy formulations of these estimates.  We also develop a Freiman-type inverse theorem for a special class of $2$-step nilpotent groups, namely the Heisenberg groups with no $2$-torsion in their centre.
\end{abstract}

\maketitle

\section{Introduction} 

The field of \emph{additive combinatorics} is concerned with the structure and size properties of sum sets such as $A+B := \{ a+b: a \in A, b \in B \}$ for various sets $A$ and $B$ (in some additive group $G$).  One also considers partial sum sets such as
$A \stackrel{E}{+} B := \{ a+b: (a,b) \in E \}$ for some\footnote{We use $E$ here instead of the more traditional $G$, as we are reserving $G$ for the ambient group.} $E \subset A \times B$.
There are many deep and important results in this theory, but we shall mention three particularly important ones.  Firstly, there are the \emph{Pl\"unnecke-Ruzsa sum-set estimates}, which roughly speaking asserts that if one sum-set such as $A+B$ is small, then other sum-sets such as $A-B$, $A+B+B$, $A+A$, etc. are also small; see e.g. \cite{plun}, \cite{ruzsa}, \cite{nathanson}, \cite{tao-vu}.  Then there is the \emph{Balog-Szemer\'edi-Gowers theorem} \cite{balog}, \cite{gowers-4}, which roughly speaking asserts that if a partial sum-set
$A \stackrel{E}{+} B$ is small for some dense subset $E$ of $A \times B$, then there are large subsets $A', B'$ of $A$, $B$
respectively whose \emph{complete} sum-set $A' + B'$ is also small.  Finally, there are \emph{inverse sum set theorems}, of which \emph{Freiman's theorem} \cite{frei} (see also \cite{bilu-freiman}, \cite{ruzsa-freiman}, \cite{chang}, \cite{nathanson}, \cite{tao-vu}) is the most famous: it asserts that if $A$ is a finite non-empty subset of a torsion-free abelian group (such as $\Z^d$) with $A+A$ small, then $A$ can be efficiently contained in the sum of $O(1)$ arithmetic progressions.
These three families of results have had many applications, perhaps most strikingly to the work of Gowers \cite{gowers-4}, \cite{gowers} on quantitative bounds for Szemer\'edi's theorem and to the work of Bourgain and co-authors \cite{bourgain-diffie}, \cite{bourgain-mordell}, \cite{bkt}, \cite{konyagin} on exponential sum estimates in finite fields.  We refer the reader to \cite{tao-vu} for a more detailed treatment of these topics.

The above results are usually phrased in the discrete setting, with $A$ and $B$ being a finite subset of an abelian group such as the
lattice $\Z^d$, and with the cardinality $|A|$ of a set $A$ used as a measure of size.  However, it is easy to transfer these discrete results to a continuous setting, for instance when $A$ and $B$ are open bounded subsets of a Euclidean space $\R^d$, 
and using Lebesgue measure $\mu(A)$ rather than cardinality to measure the size of a set.  Indeed one can pass from the continuous case to the discrete case (possibly losing constants which are exponential in the dimension $d$) by discretizing $\R^d$ to a fine lattice such as $\eps \cdot \Z^d$, applying the discrete sum-set theory, and then taking limits as $\eps \to 0$.  For similar reasons, there is little difficulty in transferring the sum-set theory to a metric entropy setting, in which the size of a set $A$
in $\R^d$ is measured using a covering number $\N_\eps(A)$. See for instance Propositions \ref{frei-cts}, \ref{frei-entropy} below for examples of these transference techniques.

In this paper we present analogues of the Pl\"unnecke-Ruzsa and Balog-Szemer\'edi theorems in the non-commutative setting, in which $G$ is now a multiplicative group and one studies the size of product sets $A \cdot B := \{ a \cdot b: a \in A, b \in B \}$ and partial product sets
$A \stackrel{E}{\cdot} B := \{ a \cdot b: (a,b) \in E \}$; the theory for inverse sumset estimates is significantly more complicated and does not seem to easily extend to the non-commutative setting (as the Fourier-analytic techniques are substantially less effective in this case), though we are able to obtain an inverse theorem for a class of Heisenberg groups, see Theorem \ref{heisen} below.  
The other two results, however, are more elementary in nature, and much of the theory carries over surprisingly easily to
this setting.  The one result which fails utterly is the Pl\"unnecke magnification inequality \cite{plun}, which is valid only for commutative graphs and has no known counterpart in the non-commutative setting.  Fortunately, one can use other, more elementary
combinatorial arguments as a substitute for the Pl\"unnecke inequalities (as was done in \cite{tao-vu}), albeit at the cost of degrading the exponents in the estimates slightly.  Another significant issue is that there is no obvious relationship between the size of $A \cdot B$ and of $B \cdot A$ in the non-commutative setting; consider for instance the case when $A$ is a subgroup of $G$, and $B$ is a right coset of $A$.  Fortunately, there is a residual relationship between the sets $A \cdot A^{-1}$ and $A^{-1} \cdot A$ in that if one product set is small then a large portion of the other product set is small (see Lemma \ref{eaa} below).  This
key observation allows us to get around the obstruction of non-commutativity and recover almost all the standard sum-set theory,
though in some cases one has to throw out a small exceptional set in order to proceed.  For instance, if $A$ is the union of a subgroup and a point, then $A \cdot A$ will be small, but higher products such as $A \cdot A \cdot A$ can be large; however, by throwing out this exceptional point one can control all products.

Finally, the passage between the discrete setting, the continuous setting, and the metric entropy setting is not as automatic
in the non-commutative setting (such as for the group $SU(2)$) as it is in the case of Euclidean spaces $\R^d$, because 
there are usually no good analogues of the discrete subgroups $\eps \cdot \Z^d$ in the general setting.  Fortunately, the
continuous and discrete theories are almost identical, so much so that we shall treat the two in a unified manner.  One can
then pass from the continuous setting to a metric entropy setting by standard volume packing arguments, provided that the metric
structure is sufficiently compatible with the group structure.  Ideally one wants the metric to be bi-invariant, but this is usually
only possible when the group $G$ is compact.  For non-compact groups such as $SL(2,\R)$, the metric entropy results that we present
here are only satisfactory when all the sets under consideration are contained in a fixed bounded set, in which case the metric structure will be approximately bi-invariant (or more precisely, the group operations are Lipschitz) and the metric balls will obey a doubling property, in which case the volume packing arguments go through without difficulty.  

Our main results are as follows.  In both the continuous and discrete setting, we classify sets of small tripling, showing that such sets are nothing more than dense subsets of a type of set that we call an \emph{approximate group}; see Theorem \ref{tripling-classify}.  As for sets of small doubling (or pairs of sets with small product set), we have a slightly different classification, showing that such sets can be covered efficiently by left or right-translates of an approximate group (Theorem \ref{energy-gleam}).  For pairs of sets with small \emph{partial} product set, we show that such sets have large \emph{intersection} with translates of an approximate group of comparable size (Theorem \ref{energy-gleam-2}).  In Section \ref{entropy-sec} we extend these results to the metric entropy setting, given some mild hypotheses on the metric.  Finally, in Section \ref{inverse-sec} we discuss the inverse product set problem (the noncommutative generalisation of the inverse sum set problem) and present a new theorem in this direction in the context of Heisenberg groups in the absence of $2$-torsion.  All of these results are \emph{polynomially reversible} in the sense that we can pass from one class of sets to an equivalent class and then back to the original class, losing only polynomial factors in the parameter $K$ (which should be thought of as a type of doubling constant).  

The author thanks Jean Bourgain for encouragement, and for raising the issue of the metric entropy case in the non-commutative setting.  He also thanks Imre Ruzsa for very detailed comments and suggestions, and Emmanuel Kowalski for corrections. This work developed from some earlier unpublished notes of the author \cite{tao-noncom}, as well as from portions of the author's book with Van Vu \cite{tao-vu}.  In particular, the discrete versions of the results here can largely be found in \cite[\S 2.7]{tao-vu}, although in some cases the proofs are assigned as exercises rather than given in full. The differences between the discrete and the continuous arguments are mostly notational in nature.

\section{Setup and notation}

We now give the unified framework in which to present the discrete and continuous non-commutative sum-set (or more precisely product-set) theory.

\begin{definition}[Multiplicative groups]  A \emph{multiplicative group} will be a topological group $G$ (thus the
group operation $(x,y) \mapsto x \cdot y$ and the inversion operation $x \mapsto x^{-1}$ are continuous), equipped with a \emph{Haar measure} $\mu$, which for us will be a non-negative Radon measure on $G$ which is invariant under left and right translation 
and inversion,  thus $\mu(x \cdot A) = \mu(A \cdot x) = \mu(A^{-1}) = \mu(A)$ for all measurable $A$ in $G$, where $x \cdot A := \{ x \cdot y: y \in G \}$, $A \cdot x := \{ y \cdot x: y \in G \}$, and $A^{-1} := \{ x^{-1}: x \in G \}$.  We denote the multiplicative identity by $1 = 1_G$.
We also make the mild non-degeneracy assumption that every non-empty open set has non-zero measure. A \emph{multiplicative set} will be any non-empty open precompact set $A$ in $G$; note that we necessarily have $0 < \mu(A) < \infty$.  
Given two multiplicative sets $A$ and $B$ we define their product $A \cdot B := \{a \cdot b: a \in A, b \in B \}$; observe that this is also a multiplicative set, as is the inverse set $A^{-1}$.
\end{definition}

\begin{remarks}
The hypotheses that a multiplicative set is open and precompact (and that $\mu$ is Radon) will allow us to avoid many 
technical issues concerning measurability and integrability, and we shall in fact not discuss these issues here.
Note that we are implicitly assuming that $G$ is locally compact, since otherwise there will be no multiplicative sets to consider.
One could weaken the translation and inversion invariance properties of the measure somewhat (so that the group operations only preserve the measure \emph{approximately}) but this would introduce a number of measure-dependent constants into the estimates below and we will not do so here.  However, such a generalisation would be useful for studying non-unimodular Lie groups.
\end{remarks}

We now give the two main examples of multiplicative groups.

\begin{example}[Discrete case]\label{discrete-ex}  Let $G$ be an abstract group (not necessarily abelian).  Then we can equip this group with the discrete topology and counting measure $\mu(A)=|A|$ to obtain a multiplicative group.  In this case, the multiplicative sets are simply the finite non-empty sets.
\end{example}

\begin{example}[Unimodular Lie group case]  Let $G$ be a finite-dimensional unimodular Lie group.  This is a finite-dimensional manifold and thus comes with a standard topology, and a standard Haar measure (defined up to a normalizing scalar).  The multiplicative sets in this case are the non-empty bounded open sets.
\end{example}

\begin{remark}
In the commutative setting one can pass between the discrete and continuous cases above by standard discretisation arguments, but the connection between
the two is less clear in the non-commutative setting.  Nevertheless we shall be able to treat both of these cases in a completely
unified manner.
\end{remark}

\begin{remark} Observe that the hypotheses on the measure $\mu$ are preserved if we multiply the measure $\mu$ by a positive constant.  Thus all the estimates we present in this paper will be invariant under this symmetry; roughly speaking, this means that the number of times $\mu$ appears on the left-hand side of an equality will always equal the number of times $\mu$ appears on the right-hand side.  (Certain quantities such as the Ruzsa distance $d(A,B)$ and the doubling constant $K$ will be dimensionless,
whereas the multiplicative energy $\Energy(A,B)$, which we define below, has the units of $\mu^3$.)
\end{remark}

Henceforth we fix the multiplicative group $G$ (and the measure $\mu$).  In the next few sections we study how 
the measure of various products such as $\mu(A \cdot B)$, $\mu(A \cdot A)$, $\mu(A \cdot A \cdot A)$, etc. of multiplicative sets are related.

We shall use the notation $X = O(Y)$, $Y = \Omega(X)$, $X \lesssim Y$ or $Y \gtrsim X$ to denote the statement 
that $X \leq CY$ for an absolute constant $C$ (not depending on the choice of group $G$ or on any other parameters).  We also use
$X \sim Y$ to denote the estimates $X \lesssim Y \lesssim X$.  If we wish to indicate dependence of the constant on
an additional parameter, we will subscript the notation appropriately, thus for instance $X \sim_n Y$ denotes that $X \leq C_n Y$ and $Y \leq C_n X$ for some $C_n$ depending on $n$.

\section{Ruzsa distance, and tripling sets}

To measure the multiplicative structure inherent in a multiplicative set $A$, or a pair $A,B$ of multiplicative
sets, it is convenient to introduce two measurements, the \emph{Ruzsa distance} and the \emph{multiplicative energy}.
In this section we focus on the Ruzsa distance and applications to sets of small tripling.

\begin{definition}[Ruzsa distance]
Let $A$ and $B$ be multiplicative sets.  We define the \emph{(left-invariant) Ruzsa distance} $d(A,B)$ to be the quantity
$$ d(A,B) := \log \frac{\mu(A \cdot B^{-1})}{\mu(A)^{1/2} \mu(B)^{1/2}}.$$
\end{definition}

We now justify the terminology ``left-invariant\footnote{One could also define a right-invariant Ruzsa distance
$\tilde d(A,B) := d(A^{-1},B^{-1}) = \log \frac{\mu(B^{-1} \cdot A)}{\mu(A)^{1/2} \mu(B)^{1/2}}$, but we will not need that notion here.} Ruzsa distance''.

\begin{lemma}[Ruzsa triangle inequality]\label{Ruzsa-triangle}  Let $A,B,C$ be multiplicative sets.  Then we have $d(A,B) \geq 0$, $d(A,B) = d(B,A)$,
and $d(A,C) \leq d(A,B) + d(B,C)$.  Also we have $d(x \cdot A, x \cdot B) = d(A,B)$ for all $x \in G$.
\end{lemma}

\begin{remark} This inequality was first established in the discrete case in \cite{ruzsa} (initially in the commutative case, but the argument extends easily to non-commutative settings).
\end{remark}

\begin{proof} From translation invariance we have $\mu(A \cdot B^{-1}) \geq \mu(A \cdot b^{-1}) =\mu(A)$ for any $b \in B$.  Similarly $\mu(A \cdot B^{-1}) \geq \mu(B^{-1}) = \mu(B)$.  Taking geometric means we obtain $\mu(A \cdot B^{-1}) \geq \mu(A)^{1/2} \mu(B)^{1/2}$ and hence $d(A,B) \geq 0$.  The symmetry property $d(A,B) = d(B,A)$ follows from the fact that $B \cdot A^{-1}$ is the
inverse of $A \cdot B^{-1}$.  Finally, to show the triangle inequality $d(A,C) \leq d(A,B) + d(B,C)$, it will suffice to show the inequality
$$ \mu(A \cdot C^{-1}) \mu(B) \leq \mu(A \cdot B^{-1}) \mu(B \cdot C^{-1}).$$
To prove this, we rewrite the right-hand side as a double integral
$$ \int_G \int_G 1_{A \cdot B^{-1}}(x) 1_{B \cdot C^{-1}}(y)\ d\mu(x) d\mu(y),$$
where $1_A$ denotes the indicator function of $A$.
Making the substitution $x = z \cdot y^{-1}$ and using the translation invariance and Fubini's theorem, we can rewrite this as
$$ \int_G \left[ \int_G 1_{A \cdot B^{-1}}(z \cdot y^{-1}) 1_{B \cdot C^{-1}}(y)\ d\mu(y)\right] d\mu(z).$$
Now if $z$ lies in $A \cdot C^{-1}$, then we have $z = a \cdot c^{-1}$ for some $a \in A$ and $c \in C$, and
then $1_{A \cdot B^{-1}}(z\cdot y^{-1}) 1_{B \cdot C^{-1}}(y) = 1$ whenever $y \in B \cdot c^{-1}$.  Since $B \cdot c^{-1}$ has the same measure as $B$, the above integral is at least as large as
$$ \int_{A \cdot C^{-1}} \mu( B )\ d\mu(z) = \mu(A \cdot C^{-1}) \mu(B)$$
and the claim follows.
\end{proof}

We caution that $d(A,A) = \log \frac{\mu(A \cdot A^{-1})}{\mu(A)}$ will usually not be zero.  Also, $d(A,B) \neq d(A^{-1}, B^{-1})$
and $d(A,B) \neq d(A \cdot x, B \cdot x)$ in general.

From the Ruzsa triangle inequality we see in particular that 
\begin{equation}\label{daa}
d(A,A) \leq 2d(A,B)
\end{equation}
for all multiplicative sets $A,B$.

For any integer $n \geq 1$ and any multiplicative set $A$, let $A^n$ denote the $n$-fold product set
$$ A^n = A \cdot \ldots \cdot A = \{ a_1 \ldots a_n: a_1,\ldots,a_n \in A \}.$$
In the commutative setting, the Pl\"unnecke-Ruzsa inequalities show that if $A^2$ is comparable in size to $A$, then $A^n$ is
also comparable in size to $A$ for any fixed $n$.  The same statement is not necessarily true in the non-commutative setting:
for instance, in the discrete setting if we take $A = H \cup \{x\}$, where $H$ is a finite subgroup of $G$ and $x$ lies outside of the normalizer of $H$, then $A^2$ has size comparable to $A$, but $A^3$ can be much larger.  However, (as was observed in \cite{Helf} in the discrete non-commutative case, although the basic argument is essentially in \cite{rt}) it turns out that once $A^3$ is under control, then so are all other combinations of $A$ and $A^{-1}$:

\begin{lemma}\label{tripling-ok} Let $A$ be a multiplicative set such that $\mu(A^3) \leq K \mu(A)$.  Then for any signs $\epsilon_1,\ldots,\epsilon_n \in \{-1,1\}$ we have
$\mu( A^{\epsilon_1} \ldots A^{\epsilon_n} ) \leq K^{O_n(1)} \mu(A)$.
\end{lemma}

\begin{proof} Let us first observe from hypothesis that $\mu(A) \leq \mu(A^2) \leq \mu(A^3) \leq K \mu(A)$.  This implies that
$d(A, A^{-1}) \leq \log K$ and $d(A^2, A^{-1}) \leq \log K$.  By the triangle inequality we thus have $d(A^2,A) = O(\log K)$, 
and thus $\mu(A \cdot A \cdot A^{-1}) \leq K^{O(1)} \mu(A)$, which implies that $d(A, A \cdot A^{-1}) = O(\log K)$.
By the triangle inequality again this implies $d(A \cdot A^{-1}, A^{-1}) = O(\log K)$, hence $\mu(A \cdot A^{-1} \cdot A)
\leq K^{O(1)} \mu(A)$.  In particular $d(A, A^{-1} \cdot A) = O( \log K)$, so again by the triangle inequality
$d(A^{-1}, A^{-1} \cdot A) = O(\log K)$, and hence $\mu(A^{-1} \cdot A \cdot A) \leq K^{O(1)} \mu(A)$.  With all these bounds
(and taking inverses) we can already establish the lemma when $n=3$, which also implies the lemma when $n < 3$.

Now we assume inductively that the lemma is already proven for all $n < n_0$ for some $n_0 \geq 4$, and wish to prove it
for $n = n_0$.  To establish the bound on $A^{\epsilon_1} \ldots A^{\epsilon_n}$ it suffices to establish the bound
$$ d( A^{\epsilon_1} \ldots A^{\epsilon_{n-2}}, A^{-\epsilon_n} \cdot A^{-\epsilon_{n-1}} ) = O_n(\log K).$$
But since the lemma is already proven for $n-1$, we have
$$ d( A^{\epsilon_1} \ldots A^{\epsilon_{n-2}}, A ) = O_n(\log K)$$
and since the lemma is already proven for $3$, we have
$$ d( A, A^{-\epsilon_n} \cdot A^{-\epsilon_{n-1}} ) = O(\log K)$$
and so the claim follows from the triangle inequality.
\end{proof}

\begin{remark} One can weaken the condition $\mu(A^3) \leq K|A|$ to $\sup_{a \in A} |A \cdot a \cdot A| \leq K |A|$; see Corollary \ref{corruz}.
\end{remark}

One can analyze the behavior of tripling sets further, by the following covering lemma.

\begin{lemma}[Ruzsa covering lemma]\label{cover}  Let $A, B$ be multiplicative sets such that $\mu(A \cdot B) \leq K \mu(A)$ (resp. $\mu(B \cdot A) \leq K \mu(A)$).  Then there exists a finite set $X$ contained inside $B$ of cardinality at most $K$ such that $B \subseteq A^{-1} \cdot A \cdot X$
(resp. $B \subseteq X \cdot A \cdot A^{-1}$).
\end{lemma}

\begin{remark} For the commutative version of this lemma (in discrete or continuous settings), see \cite{ruzsa-group}, \cite{milman}, \cite{tao-vu}.
\end{remark}

\begin{proof}  By the reflection invariance of $\mu$ it will suffice to prove the claim when $\mu(A \cdot B) \leq K \mu(A)$.
Let $X$ be a subset of $B$ with the property that the sets $A \cdot x$ for $x \in X$ are disjoint.  Since $\mu(A \cdot B) \leq K \mu(A)$ we see that such a set $X$ must have cardinality at most $K$. Now let $X$ be such a set which is maximal with respect to set inclusion (which one can construct for instance using the greedy algorithm).  Then for any $b \in B$ we must have $A \cdot b$ intersecting $A \cdot x$ for some $x \in X$, which
implies that $b \in A^{-1} \cdot A \cdot X$.  The claim follows.
\end{proof}

We can now give a classification of sets of small tripling.

\begin{definition}[Approximate groups]\label{agdef} A multiplicative set $H$ is said to be a \emph{$K$-approximate group} if
it is symmetric (so $H^{-1} = H$) and there exists a finite symmetric set $X \subset H^2$ of cardinality at most $K$ such
that $H \cdot H \subseteq X \cdot H$.
\end{definition}

\begin{remark}  In \cite{tao-vu}, the additional condition $X \cdot H \subset H \cdot X \cdot X$ was also imposed.  Our methods for constructing approximate groups also yield this additional property (see Theorem \ref{tripling-classify} below), though we will not use this additional hypothesis in our arguments and thus omit it from the definition of an approximate group.
\end{remark}

Note from the symmetry assumptions that if $H \cdot H \subset X \cdot H$, then $H \cdot H \subset H \cdot X$ also.
Iterating this we see that $H^n \subseteq X^{n-1} \cdot H, H \cdot X^{n-1}$ for all $n \geq 1$.  Thus approximate groups have small
tripling.  It turns out that this is essentially the only way that a multiplicative set can have small tripling:

\begin{theorem}\label{tripling-classify}
  Let $K \geq 1$, and let $A$ be a multiplicative set.  Then the following three statements are equivalent, in the sense that if one of them holds for one choice of implied constant in the $O()$ and $\lesssim$ notation, then the other statements hold for a different choice of implied constant in the $O()$ and $\lesssim$ notation:
\begin{itemize}
\item[(i)] We have the tripling bound $\mu(A^3) \lesssim K^{O(1)} \mu(A)$.
\item[(ii)] We have $\mu( A^{\epsilon_1} \ldots A^{\epsilon_n} ) \sim_n K^{O_n(1)} \mu(A)$ for all $n \geq 1$ and all signs $\epsilon_1, \ldots, \epsilon_n \in \{-1,+1\}$.
\item[(iii)] There exists a $O(K^{O(1)})$-approximate group $H$ of size $\mu(H) \sim K^{O(1)} \mu(A)$ 
which contains $A$.
\end{itemize}
\end{theorem}

\begin{proof}  The implication (i) $\implies$ (ii) is Lemma \ref{tripling-ok}, while the reverse implication (ii) $\implies$ (i) is trivial.
The implication (iii) $\implies$ (i) is also trivial:
$$ \mu(A^3) \leq \mu(H^3) \lesssim K^{O(1)} \mu(H) \lesssim K^{O(1)} \mu(A).$$
It remains to show that (i) implies (iii).  Set $H_0 := A \cup \{1\} \cup A^{-1}$ and $H := H_0^3$, then from Lemma \ref{tripling-ok} we see that $\mu(H) \sim K^{O(1)} \mu(A)$.  Clearly $H$ also contains $A$, so it
remains to show that $H$ is a $O(K^{O(1)})$-approximate group.  Certainly $H$ is symmetric.  From Lemma \ref{tripling-ok} we have
$\mu( H_0 \cdot H^2 ) \lesssim K^{O(1)} \mu(A)$, and hence from Lemma \ref{cover} we can find a finite set $Y$ in $H^2$
of cardinality $O(K^{O(1)})$ such that 
$$H^2 \subseteq H_0^{-1} \cdot H_0 \cdot Y \subseteq H \cdot Y.$$
If we set $X := Y \cup Y^{-1}$ then $X$ is symmetric and (from symmetry of $H$) we conclude that $H^2 \subseteq H \cdot X, X \cdot H$.  But since $X$ is contained in $H^2$, we also have
$$ H \cdot X \subseteq H \cdot H^2 = H^2 \cdot H \subseteq X \cdot H \cdot H \subseteq X \cdot X \cdot H$$
and similarly $H \cdot X \subseteq H \cdot X \cdot X$.  The claim follows.
\end{proof}

An inspection of the above proof reveals the following more precise implication of (i) from (iii):

\begin{corollary}\label{tripling-better} Let $A$ be a multiplicative set such that $\mu(A^3) \leq K \mu(A)$.  Then
the set $H := (A \cup \{1\} \cup A^{-1})^3$ is a $O(K^{O(1)})$-approximate group.  In particular, if $A$ is symmetric and contains $1$, then $A^3$ is a $O(K^{O(1)})$-approximate group.
\end{corollary}

\section{Convolution and multiplicative energy}

To study sets of small doubling, rather than small tripling, it is convenient to introduce another measure of multiplicative
structure between two multiplicative sets, namely the \emph{multiplicative energy}.  Given two absolutely integrable functions $f, g$ on $G$, we define their \emph{convolution}
$f*g$ in the usual manner as
$$ f*g(x) := \int_G f(y) g(y^{-1} x)\ d\mu(y) = \int_G f(x y^{-1}) g(y)\ d\mu(y).$$
As is well-known, convolution is bilinear and associative (though not necessarily commutative), and the convolution of 
two absolutely integrable functions is continuous.  Convolution is not commutative in general, but we do have the identity
\begin{equation}\label{fg0}
f*g(1) = g*f(1),
\end{equation}
which reflects the fact that $xx^{-1} = x^{-1} x = 1$.
If $f$ is supported on $A$ and $g$ is supported on $B$ then $f*g$ is supported on
$A \cdot B$.  If we use $\tilde f(x) := f(x^{-1})$ to denote the reflection of $f$, we observe the reflection property
\begin{equation}\label{reflect}
\widetilde{f*g} = \tilde g * \tilde f
\end{equation}
(which reflects the fact that $(x \cdot y)^{-1} = y^{-1} \cdot x^{-1}$ for $x,y \in G$)
and the trace formula
\begin{equation}\label{fg}
\int_G f(x) g(x)\ d\mu(x) = f * \tilde g(1) = \tilde g * f(1) = \tilde f * g(1) = f * \tilde g(1).
\end{equation}
If $A$ is a multiplicative set, we use $1_A$ to denote the indicator function of $A$; note that this is an absolutely integrable function.

\begin{definition}  Let $A, B$ be multiplicative sets.  We define the \emph{multiplicative energy} $\Energy(A,B)$ between these two sets
to be the quantity
$$ \Energy(A,B) := \int_G [1_A * 1_B(x)]^2\ d\mu(x).$$
\end{definition}

\begin{remark} In the discrete setting (Example \ref{discrete-ex}), we have
\begin{equation}\label{energy-discrete}
\Energy(A,B) = |\{ (a,b,a',b') \in A \times B \times A \times B: a \cdot b = a' \cdot b' \}|.
\end{equation}
In the notation of Gowers \cite{gowers-4}, the quantity $\Energy(A,B)$ thus counts the number of \emph{multiplicative quadruples} in
$A \times B \times A \times B$.  In the commutative case, this quantity has a useful representation in terms of the Fourier transform
(see e.g. \cite{gowers-4}, \cite{tao-vu}), which has a number of applications, for instance in proving Freiman's inverse sumset theorem.  In the non-commutative case it is also possible to use the non-commutative Fourier transform to represent this energy, but the resulting formulae are not as tractable.  In particular, no analogue of Freiman's theorem is currently known in the general non-commutative setting.  Fortunately, we will be able to use the properties of the convolution algebra (most notably its associativity, the reflection property \eqref{reflect}, and the trace property \eqref{fg}) to compensate for the lack of a convenient Fourier-analytic description of the energy; in particular, we will not need to understand the representation theory of the underlying group $G$.
\end{remark}

A simple application of Fubini's theorem and change of variables shows that
\begin{equation}\label{fubini}
 \int_G 1_A * 1_B(x)\ d\mu(x) = \mu(A) \mu(B)
\end{equation}
and
\begin{equation}\label{1ab}
 1_A * 1_B(x) \leq \min( \mu(A), \mu(B) ) \leq \mu(A)^{1/2} \mu(B)^{1/2}
 \end{equation}
and hence by H\"older's inequality we have the upper bound
\begin{equation}\label{eab} \Energy(A,B) \leq \mu(A)^{3/2} \mu(B)^{3/2}.
\end{equation}
Also, since $1_A * 1_B$ is supported on $A \cdot B$, we have a lower bound from \eqref{fubini} and Cauchy-Schwarz:
\begin{equation}\label{eab-lower}
 \Energy(A,B) \geq \frac{\mu(A)^2 \mu(B)^2}{\mu(A \cdot B)}.
\end{equation}
In general, we do not have $\Energy(A,B) = \Energy(B,A)$ (although one can use \eqref{reflect} to show that $\Energy(A,B) = \Energy(B^{-1},A^{-1})$).
On the other hand, we do have the following important identity, which can be viewed as a weak form of commutativity.

\begin{lemma}\label{eaa} For any multiplicative set $A$, we have $\Energy(A,A^{-1}) = \Energy(A^{-1},A)$.
\end{lemma}

\begin{remark}\label{aai} This identity is especially striking since there is no relation between the size of $A \cdot A^{-1}$ and $A^{-1} \cdot A$ in general.  For instance, if $H$ is a multiplicative set which is also a subgroup of $G$, and $A := (x \cdot H) \cup H$ for some $x$ not in the normaliser of $H$, then $A \cdot A^{-1}$ has about the same size as $H$, but $A^{-1} \cdot A$ can be much larger.  In the discrete case (Example \ref{discrete-ex}) one can prove this lemma using the identity \eqref{energy-discrete} and the observation that $a \cdot b = a' \cdot b'$ if and only if $b \cdot (b')^{-1} = a^{-1} \cdot a'$.
\end{remark}

\begin{proof}  From \eqref{fg}, \eqref{reflect} (and associativity) we have
$$ \Energy(A,A^{-1}) = 1_A * 1_{A^{-1}} * \widetilde{1_A * 1_{A^{-1}}}(1) = 1_A * 1_{A^{-1}} * 1_A * 1_{A^{-1}}(1).$$
Similarly we have $\Energy(A^{-1},A) = 1_{A^{-1}} * 1_A * 1_{A^{-1}} * 1_A(1)$.  The claim then follows from \eqref{fg0}.
\end{proof}

This lemma has the following useful consequence.

\begin{proposition}\label{musprop}  Let $A$ be a multiplicative set such that $\mu(A \cdot A^{-1}) \leq K \mu(A)$.
Then there exists a symmetric multiplicative set $S$ such that $\mu(S) \geq \mu(A)/2K$ and
\begin{equation}\label{mus}
\mu( A \cdot S^n \cdot A^{-1}) \leq 2^n K^{2n+1} \mu(A)
\end{equation}
for all integers $n \geq 1$.
\end{proposition}

\begin{proof}  From Lemma \ref{eaa} and \eqref{eab-lower} we have
\begin{align*}
\int_G \mu( A \cap (A \cdot x) )^2\ d\mu(x) &= 1_{A^{-1}} * 1_A(x)^2\ d\mu(x)\\
&= \Energy(A^{-1},A) \\
&= \Energy(A,A^{-1}) \\
&\geq \mu(A)^4 / \mu(A \cdot A^{-1}) \\
&\geq \mu(A)^3/K.
\end{align*}
Now we define $S$ as
$$ S := \{ x \in G: \mu( A \cap (A \cdot x) ) > \mu(A)/2K \}.$$
It is easy to see that $S$ is a symmetric multiplicative set.
From \eqref{fubini} we see that
$$ \int_G \mu(A \cap (A \cdot x))\ d\mu(x) = \mu(A)^2$$
and thus
$$ \int_{G \backslash S} \mu(A \cap (A \cdot x))^2\ d\mu(x) \leq \mu(A)^3/2K.$$
Subtracting this from the preceding estimate, we conclude
$$ \int_S \mu(A \cap (A \cdot x))^2\ d\mu(x) \geq \mu(A)^3/2K.$$
Bounding $\mu(A \cap (A \cdot x))$ by $\mu(A)$, we conclude in particular that $\mu(S) \geq \mu(A)/2K$.

It remains to prove \eqref{mus}.  Let us consider the quantity
\begin{equation}\label{amu}
 \int_{(A \cdot A^{-1})^{n+1}} 1_{A \cdot S^n \cdot A^{-1}}(y_0 \ldots y_n)\ d\mu(y_0) \ldots d\mu(y_n).
\end{equation}
On the one hand, this quantity is clearly bounded above by $\mu(A \cdot A^{-1})^{n+1} \leq K^{n+1} \mu(A)$.  Now
let us obtain a lower bound.  We rewrite this quantity as
$$ \int_{A \cdot S^n \cdot A^{-1}} \int_{(A \cdot A^{-1})^n} 1_{A \cdot A^{-1}}(y_{n-1}^{-1} \ldots y_0^{-1} x)\ 
d\mu(y_0) \ldots d\mu(y_{n-1}) d\mu(x).$$
Suppose that we can show that
\begin{equation}\label{xyx}
\int_{(A \cdot A^{-1})^n} 1_{A \cdot A^{-1}}(y_{n-1}^{-1} \ldots y_0^{-1} x)\ 
d\mu(y_0) \ldots d\mu(y_{n-1}) \geq (\mu(A)/2K)^n
\end{equation}
for all $x \in A \cdot S^n \cdot A^{-1}$  Then we can bound \eqref{amu} from below by
$\mu(A \cdot S^n \cdot A^{-1}) (\mu(A)/2K)^n$, which will establish \eqref{mus}.

It remains to show \eqref{xyx}.  Let $x \in A \cdot S^n \cdot A^{-1}$ be arbitrary.  We can write
$x = a_0 s_1 \ldots s_n b_{n+1}^{-1}$ where $a_0, b_{n+1} \in A$ and $s_1,\ldots,s_n \in S$.  If we make the successive change of variables
$$ y_0 = a_0 b_1^{-1}; \quad y_1 = b_1 s_1 b_2^{-1}; \quad \ldots \quad; y_{n-1} = b_{n-1} s_{n-1} b_n^{-1}$$
then we observe that
$$ y_{n-1}^{-1} \ldots y_0^{-1} x = b_n s_n b_{n+1}^{-1}$$
and we can rewrite the left-hand side of \eqref{xyx} as
$$\int_{G^n} 1_{A \cdot A^{-1}}( a_0 b_1^{-1} ) \prod_{i=1}^{n} 1_{A \cdot A^{-1}}(b_i s_i b_{i+1}^{-1})\ d\mu(b_1) \ldots d\mu(b_n).$$
Note that if $b_1,\ldots,b_n \in A$ and $b_1 s_1, \ldots, b_n s_n \in A$ then the integrand here is equal to $1$.  Hence we can
bound \eqref{xyx} from below by
$$ \int_{G^n} \prod_{i=1}^n 1_A(b_i) 1_A(b_i s_i)\ d\mu(b_1) \ldots d\mu(b_n)
= \prod_{i=1}^n \mu( A \cap (A \cdot s_i) ).$$
The claim now follows from the definition of $S$.
\end{proof}

Now we can classify sets of small doubling using approximate groups.

\begin{theorem}\label{energy-gleam}  Let $K \geq 1$, and let $A,B$ be multiplicative sets.  Then the following two statements are equivalent, in the sense that if one of them holds for one choice of implied constant in the $O()$ and $\lesssim$ notation, then the other statement holds for a different choice of implied constant in the $O()$ and $\lesssim$ notation:
\begin{itemize}
\item[(i)] We have the product bound $\mu(A \cdot B) \lesssim K^{O(1)} \mu(A)^{1/2} \mu(B)^{1/2}$ (or equivalently, $d(A,B^{-1}) \lesssim 1 + \log K$).
\item[(ii)] There exists a $O(K^{O(1)})$-approximate group $H$ of size $\mu(H) \lesssim K^{O(1)} \mu(A)^{1/2} \mu(B)^{1/2}$ and a finite set $X$ of cardinality $O(K^{O(1)})$ such that $A \subset X \cdot H$ and $B \subset H \cdot X$.
\end{itemize}
\end{theorem}

\begin{proof}
The implication (ii) $\implies$ (i) is trivial:
$$ \mu(A \cdot B) \leq \mu( X \cdot H \cdot H \cdot X ) \leq |X|^2 \mu(H^2) \lesssim K^{O(1)}  \mu( H ) \lesssim K^{O(1)} \mu(A)^{1/2} \mu(B)^{1/2}.$$
Now we show that (i) implies (ii).  From \eqref{daa} we have $d(A,A) \lesssim 1 + \log K$, thus $\mu(A \cdot A^{-1}) \lesssim K^{O(1)}  \mu(A)$.  Applying Proposition \ref{musprop}, we obtain a symmetric multiplicative set $S$ with $\mu(S) \gtrsim K^{O(1)} \mu(A)$ such that
\begin{equation}\label{muscle}
\mu( A \cdot S^3 \cdot A^{-1}) \lesssim K^{O(1)} \mu(A).
\end{equation}
In particular we have $\mu(S), \mu(A \cdot S) \lesssim K^{O(1)} \mu(A)$, which implies that $d(A,S) \lesssim 1 + \log K$.  We
also see that $\mu(S^3) \lesssim K^{O(1)} \mu(S)$, so by Theorem \ref{tripling-classify} we can find a $O(K^{O(1)})$-approximate group $H$ of size $O(K^{O(1)}) \mu(A)$ which contains $S$.  In particular we have $\mu(S \cdot H) \leq \mu(H^2) \lesssim K^{O(1)} \mu(A)$, and hence $d(S,H) \lesssim 1 + \log K$.  From the triangle inequality we conclude $d(A,H) \lesssim 1 + \log K$, thus
$\mu(A \cdot H) \lesssim K^{O(1)} \mu(A)$.  Applying Lemma \ref{cover}, there exists a finite set $Y$ of cardinality $O(K^{O(1)})$
such that $A \subset Y \cdot H \cdot H$; since $H$ is a $O(K^{O(1)})$-approximate group, we can thus find another finite set $Z$
of cardinality $O(K^{O(1)})$ such that $A \subset Z \cdot H$.  Now since $d(A,B^{-1}) \lesssim 1 + \log K$, the triangle inequality
also gives $d(B^{-1}, H) \lesssim 1 + \log K$, so by arguing as before we can find a finite set $W$ of cardinality $O(K^{O(1)})$ such
that $B^{-1} \subseteq W \cdot H$.  The claim now follows by taking $X := Z \cup W^{-1}$.
\end{proof}

One can of course specialize this to theorem to the case $A=B$, to characterize sets of small doubling:

\begin{corollary}\label{energy-cor}  Let $K \geq 1$, and let $A$ be a multiplicative set.  Then the following two statements are equivalent, in the sense that if one of them holds for one choice of implied constant in the $O()$ and $\lesssim$ notation, then the other statement holds for a different choice of implied constant in the $O()$ and $\lesssim$ notation:
\begin{itemize}
\item[(i)] We have the product bound $\mu(A \cdot A) \lesssim K^{O(1)} \mu(A)$ (or equivalently, $d(A,A^{-1}) \lesssim 1 + \log K$).
\item[(ii)] There exists a $O(K^{O(1)})$-approximate group $H$ of size $\mu(H) \lesssim K^{O(1)} \mu(A)$ and a finite set $X$ of cardinality $O(K^{O(1)})$ such that $A \subset (X \cdot H) \cap (H \cdot X)$.
\end{itemize}
\end{corollary}

As one consequence of this corollary, we can obtain the following strengthening of Lemma \ref{tripling-ok} which was conjectured to us by Imre Ruzsa (private communication):

\begin{corollary}\label{corruz}  Let $A$ be a multiplicative set such that $\mu(A \cdot a \cdot A) \leq K\mu(A)$ for all $a \in A$, and such that $\mu(A^2) \leq K\mu(A)$.  Then $\mu(A^3) \lesssim K^{O(1)} \mu(A)$, and in particular the conclusions of Lemma \ref{tripling-ok} hold (with slightly worse implied constants).
\end{corollary}

\begin{proof} By Corollary \ref{energy-cor}, we may find a $O(K^{O(1)})$-approximate group $H$ with $\mu(H) \lesssim K^{O(1)} \mu(A)$ and a finite set $X$ of cardinality $O(K^{O(1)})$ such that $A \subset X \cdot H, H \cdot X$.  By removing useless elements of $X$ if necessary we may assume that $X \subset (A \cdot H) \cup (H \cdot A)$.
Then
$$ \mu(A^3) \leq \mu( X \cdot H \cdot H \cdot X \cdot H ) \lesssim K^{O(1)} \mu( H \cdot H \cdot X \cdot H ).$$
But by Definition \ref{agdef}, $H \cdot H$ is covered by $O(K^{O(1)})$ left-translates of $H$, thus
$$ \mu(A^3) \lesssim K^{O(1)} \mu(H \cdot X \cdot H) \lesssim K^{O(1)} \sup_{x \in (A \cdot H) \cup (H \cdot A)}
\mu( H \cdot x \cdot H)$$
where the last inequality follows from the properties of $X$.  Thus it suffices to show that
$$ \mu( H \cdot x \cdot H) \lesssim K^{O(1)} \mu(A)$$
for all $x \in A \cdot H$ or $x \cdot H \cdot A$.  Splitting $x$ into factors in $H$ and $A$ and noting once again that $H \cdot H$ can be covered by $O(K^{O(1)})$ left-translates (and hence right-translates, by symmetry) of $H$, we reduce to showing that
\begin{equation}\label{mah}
\mu( H \cdot a \cdot H) \lesssim K^{O(1)} \mu(A)
\end{equation}
for all $a \in A$.

Fix $a$.  We already know that $\mu(A \cdot a \cdot A) \leq K \mu(A)$, thus
$$ d( A, A^{-1} \cdot a^{-1} ) \leq \log K.$$
On the other hand, since 
$$ \mu( A \cdot H ) \leq \mu( X \cdot H \cdot H ) \lesssim K^{O(1)} \mu(H) $$
and $\mu(H) \lesssim K^{O(1)} \mu(A)$ we see that
$$ d( A, H ) \lesssim 1+\log K.$$
By the triangle inequality we thus have
$$ d( H, A^{-1} \cdot a^{-1} ) \lesssim 1 + \log K$$
or equivalently
$$ d( H \cdot a, A^{-1} ) \lesssim 1 + \log K.$$
Now
$$ \mu( H \cdot A ) \leq \mu( H \cdot H \cdot X ) \lesssim K^{O(1)} \mu(H) $$
so by arguing as before we have
$$ d( H, A^{-1}) \lesssim 1+\log K$$
and thus by the triangle inequality
$$ d( H \cdot a, H ) \lesssim 1 + \log K$$
and the claim \eqref{mah} follows.
\end{proof}

\section{The Balog-Szemer\'edi-Gowers theorem}

In this section we develop with the non-commutative, continuous analogue of Balog-Szemer\'edi-Gowers theory.
We first give a preliminary version of this lemma, in which we start with $1/K$ of a product set $A \cdot B$ being under control, and
end up with $1-\eps$ of another product set $(A') \cdot (A')^{-1}$ being under control.  

\begin{lemma}[Weak Balog-Szemer\'edi-Gowers theorem]\label{weak-bsg}  Let $A, B, C$ be multiplicative sets such that
$$ \mu(C) \leq K' \mu(A)^{1/2} \mu(B)^{1/2}$$
and
$$ \mu \otimes \mu( \{ (a,b) \in A \times B: a \cdot b \in C \} ) \geq \mu(A) \mu(B) / K$$
for some $K, K' \geq 1$, where $\mu \otimes \mu$ denotes product measure on $G \times G$.  Let $0 < \eps < 1$.
Then there exists a multiplicative set $A'$ contained in $A$, and a multiplicative set $D$, such that
\begin{equation}\label{muap}
 \mu(A') \geq \frac{\mu(A)}{\sqrt{2} K}
\end{equation}
and
\begin{equation}\label{mud}
\mu(D) \leq \frac{2(KK')^2}{\eps} \mu(A)
\end{equation}
and
\begin{equation}\label{mumua}
\mu \otimes \mu( \{ (a,a') \in A' \times A': a \cdot (a')^{-1} \in D \} \geq (1-\eps) \mu(A')^2.
\end{equation}
\end{lemma}

\begin{proof}  By hypothesis on $C$ we have
$$ \int_B (\int_A 1_C(a \cdot b)\ d\mu(a)) d\mu(b) \geq \mu(A) \mu(B) / K.$$
By Cauchy-Schwarz we conclude that
$$ \int_B (\int_A 1_C(a \cdot b)\ d\mu(a))^2 d\mu(b) \geq \mu(A)^2 \mu(B) / K^2$$
which we rearrange as
$$ \int_A \int_A (\int_B 1_C(a \cdot b) 1_C(a' \cdot b)\ d\mu(b))\ d\mu(a) d\mu(a') \geq \mu(A)^2 \mu(B) / K^2.$$
Let $\Omega \subset A \times A$ be the set of all $(a,a')$ such that
$$ \int_B 1_C(a \cdot b) 1_C(a' \cdot b)\ d\mu(b) \leq \frac{\eps}{2K^2} \mu(B)$$
then we clearly have
$$ \int_A \int_A 1_\Omega(a,a') (\int_B 1_C(a \cdot b) 1_C(a' \cdot b)\ d\mu(b))\ d\mu(a) d\mu(a') \leq \eps \mu(A)^2 \mu(B) / 2K^2$$
and hence
$$ \int_A \int_A (1 - \frac{1}{\eps} 1_\Omega(a,a')) 
(\int_B 1_C(a \cdot b) 1_C(a' \cdot b)\ d\mu(b))\ d\mu(a) d\mu(a') \geq \mu(A)^2 \mu(B) / 2K^2.$$
We rewrite this as
$$ \int_B (\int_A \int_A (1 - \frac{1}{\eps} 1_\Omega(a,a')) 1_C(a \cdot b) 1_C(a' \cdot b)\ d\mu(a) d\mu(a')) d\mu(b) 
\geq \mu(A)^2 \mu(B) / 2K^2$$
and hence by the pigeonhole principle there exists $b \in B$ such that
$$ \int_A \int_A (1 - \frac{1}{\eps} 1_\Omega(a,a')) 1_C(a \cdot b) 1_C(a' \cdot b)\ d\mu(a) d\mu(a')
\geq \mu(A)^2 / 2K^2.$$
If we fix this $b$ and set $A' := \{ a \in A: (a,b) \in E\}$, we conclude that
$$ \mu(A')^2 \geq \int_{A'} \int_{A'} (1 - \frac{1}{\eps} 1_\Omega(a,a'))\ d\mu(a) d\mu(a')
\geq \mu(A)^2 / 2K^2 \geq 0$$
which in particular implies \eqref{muap}.  Also we see from the above inequality that
$$ \mu \otimes \mu( (A' \times A') \cap \Omega ) \leq \eps \mu(A')^2.$$
Thus if we define
$$ D := \{ a \cdot (a')^{-1}: a, a' \in A'; (a,a') \not \in \Omega \}$$
then we have \eqref{mumua}.  Now suppose that $d = a \cdot (a')^{-1}$ lies in $D$ for some
$a, a' \in A'$ and $(a,a') \not \in \Omega$.  From definition of $\Omega$ we have
$$ \int_B 1_C(a \cdot b) 1_C(a' \cdot b)\ d\mu(b) > \frac{\eps}{2K^2} \mu(B)$$
and hence by the substitution $c := a' \cdot b$
$$ \int_G 1_C(d \cdot c) 1_C(c)\ d\mu(c) > \frac{\eps}{2K^2} \mu(B).$$
Integrating this over all $d \in D$, we obtain
$$ \int_G \int_G 1_C(d \cdot c) 1_C(c)\ d\mu(c) d\mu(d) > \frac{\eps}{2K^2} \mu(B) \mu(D).$$
Using Fubini's theorem and making the change of variables $c' = d \cdot c$ we see that the left-hand side
is just $\mu(C)^2 \leq (K')^2 \mu(A) \mu(B)$, and \eqref{mud} follows.
\end{proof}

Now we extend the Balog-Szemer\'edi-Gowers theorem to pairs $A,B$ of multiplicative sets (of comparable size).

\begin{theorem}[Balog-Szemer\'edi-Gowers theorem]\label{bsg} Let $A, B$ be multiplicative sets such that $\Energy(A,B) \geq \mu(A)^{3/2} \mu(B)^{3/2} / K$.  Then there exist multiplicative sets $A''', B'''$ contained in $A$, $B$ respectively such that
$\mu(A''') \geq \frac{\mu(A)}{8 \sqrt{2} K}$, $\mu(B''') \geq \frac{\mu(B)}{8K}$, and 
$\mu(A''' \cdot B''') \lesssim K^8 \mu(A)^{1/2} \mu(B)^{1/2}$.
\end{theorem}

\begin{proof}  By hypothesis we have
$$ \int_G (1_A * 1_B(x))^2\ d\mu(x) \geq \mu(A)^{3/2} \mu(B)^{3/2} / K.$$
If we let $C$ denote the (open, precompact) set
$$ C := \{ x \in G: 1_A * 1_B(x) > \mu(A)^{1/2} \mu(B)^{1/2} / 2K\}$$
then we see from \eqref{fubini} that
$$ \int_{G \backslash C} (1_A * 1_B(x))^2\ d\mu(x) \leq \mu(A)^{3/2} \mu(B)^{3/2} / 2K$$
and hence
\begin{equation}\label{intc}
\int_C (1_A * 1_B(x))^2\ d\mu(x) \geq \mu(A)^{3/2} \mu(B)^{3/2} / 2K.
\end{equation}
In particular, $C$ is non-empty (and is thus a multiplicative set), while from \eqref{fubini} and Markov's inequality we have
\begin{equation}\label{muc}
\mu(C) \leq 2K \mu(A)^{1/2} \mu(B)^{1/2}.
\end{equation}
Also, from \eqref{intc} and \eqref{1ab} we have
$$ \int_C 1_A * 1_B(x)\ d\mu(x) \geq \mu(A) \mu(B) / 2K.$$
By Fubini's theorem and a change of variables, the left-hand side can be rearranged as
$$ \int_A (\int_B 1_C(ab)\ d\mu(b)) d\mu(a).$$
If we thus let $A'$ be the (open precompact) subset of $A$ defined by
$$ A' := \{ a \in A: \int_B 1_C(ab)\ d\mu(b) > \mu(B) / 4K \}$$
then 
$$ \int_{A \backslash A'} (\int_B 1_C(ab)\ d\mu(b)) d\mu(a) \leq \mu(A) \mu(B) / 4K$$
and hence
\begin{equation}\label{apb}
\int_{A'} (\int_B 1_C(ab)\ d\mu(b)) d\mu(a) \geq \mu(A) \mu(B) / 4K.
\end{equation}
In particular, $A'$ is non-empty (and is thus a multiplicative set).  Using the trivial bound
$\int_B 1_C(ab)\ d\mu(b) \leq \mu(B)$, we also see that
$$ \mu(A') \geq \mu(A) / 4K.$$
Let us thus write $\mu(A) = L \mu(A')$ for some $1 \leq L \leq 4K$.  From \eqref{muc} we then have
$$
\mu(C) \leq 2K L^{1/2} \mu(A')^{1/2} \mu(B)^{1/2}
$$
while from \eqref{apb} we have
$$
\mu \otimes \mu( \{ (a,b) \in A' \times B: a \cdot b \in C \} )
\geq \mu(A') \mu(B) L / 4K.
$$
We can thus apply Lemma \ref{weak-bsg} (with $\eps := 1/32K$, and with $A, K, K'$ replaced by $A', 4K/L, 2KL^{1/2}$) to find a multiplicative set $A''$ contained in $A'$ (and hence in $A$), and a multiplicative set $D$, such that
$$ \mu(A'') \geq \frac{\mu(A') L}{4\sqrt{2} K} = \frac{\mu(A)}{4 \sqrt{2} K}$$
and
$$ \mu(D) \leq \frac{2(4K/L)^2(2KL^{1/2})^2}{1/32K} \mu(A') \lesssim K^5 \mu(A') / L.$$
In particular, we have
\begin{equation}\label{mudd}
 \mu(D) \lesssim K^6 \mu(A'') / L^2 \lesssim K^6 \mu(A'').
\end{equation}
Also we have
$$
\mu \otimes \mu( \{ (a,a') \in A'' \times A'': a \cdot (a')^{-1} \in D \} \geq (1-\frac{1}{32K}) \mu(A'')^2.
$$
We can rewrite the latter estimate as
$$ \int_{A''} \mu( \{ a' \in A'': a \cdot (a')^{-1} \not \in D \} )\ d\mu(a) \leq \frac{1}{32K} \mu(A'')^2$$
so if we set
$$ A''' := \{ a \in A'': \mu( \{ a' \in A'': a \cdot (a')^{-1} \not \in D \} ) \leq \frac{1}{16K} \mu(A'') \}$$
then by Markov's inequality we have
$$ \mu(A''') \geq \mu(A'')/2 \geq \frac{\mu(A)}{2^{3.5} K}.$$
Since $A''$ is a subset of $A'$, we have
$$ \int_B 1_C(ab)\ d\mu(b) > \mu(B) / 4K  \hbox{ for all } a \in A''$$
and hence upon integrating in $a$ and Fubini's theorem
$$ \int_B (\int_{A''} 1_C(ab)\ d\mu(a)) d\mu(b) \geq \mu(A'') \mu(B) / 4K.$$
Hence if we define the (open precompact) subset $B'''$ of $B$ by
$$ B''' := \{ b \in B: \int_{A''} 1_C(ab)\ d\mu(a) > \mu(A'') / 8K \}$$
then we have by similar arguments to before that
$$ \int_{B'''} (\int_{A''} 1_C(ab)\ d\mu(a)) d\mu(b) \geq \mu(A'') \mu(B) / 8K;$$
since $\int_{A''} 1_C(ab)\ d\mu(a) \leq \mu(A'')$, we have in particular that
$$ \mu(B''') \geq \mu(B) / 8K.$$
In particular $B'''$ is non-empty and is hence a multiplicative set.  

Now let $c = a b$ for some $a \in A'''$ and $b \in B'''$.  From definition of $B'''$ we have
$$ \mu( \{ a' \in A'': a'b \in C \} ) > \mu(A'')/8K$$
while from definition of $A'''$ we have
$$ \mu( \{ a' \in A'': a \cdot (a')^{-1} \not \in D \} ) \leq \frac{1}{16K} \mu(A'').$$
Thus
$$ \mu( \{ a' \in A'': a'b \in C, a \cdot (a')^{-1} \in D \} ) > \mu(A'')/16K$$
so by setting $x := a' b$ (so that $a \cdot (a')^{-1} = cx^{-1}$) we have
$$ \int_G 1_C(x) 1_D(cx^{-1})\ d\mu(x) > \mu(A'')/ 16 K.$$
Integrating this over all $c \in A''' \cdot B'''$ we conclude
$$ \int_G \int_G 1_C(x) 1_D(cx^{-1})\ d\mu(x) d\mu(c) \geq \mu(A'') \mu(A''' \cdot B''')/ 16 K.$$
But the left-hand side is $\mu(C) \mu(D)$, so we have
$$ \mu(A''' \cdot B''') \lesssim \frac{K \mu(C) \mu(D)}{\mu(A'')}.$$
Applying \eqref{muc}, \eqref{mudd}, we conclude
$$ \mu(A''' \cdot B''') \lesssim K^8 \mu(A)^{1/2} \mu(B)^{1/2} $$
as desired.
\end{proof}

\begin{remark} There are a number of variants of this theorem, for instance one could replace the hypothesis that $\Energy(A,B)$ is large
by the hypothesis that a partial product set $A \stackrel{E}{\cdot} B$ is small (with a suitable largeness hypothesis on $E$).
One can then refine the above theorem by requiring the additional conclusion that $E$ has large intersection with $A' \times B'$; see for instance \cite{lruzsa}, \cite{borg:high-dim} for some examples of this type of refinement. Other variants of
the lemma and its proof can be found in \cite{ssv}, \cite{chang-er}.  The power of $7$ can probably be lowered further but we shall not attempt to do so here.
\end{remark}

This gives us a characterisation of pairs of multiplicative sets of large multiplicative energy.

\begin{theorem}\label{energy-gleam-2}  Let $K \geq 1$, and let $A,B$ be multiplicative sets.  Then the following four statements are equivalent, in the sense that if one of them holds for one choice of implied constant in the $O()$ and $\lesssim$ notation, then the other statements hold for a different choice of implied constant in the $O()$ and $\lesssim$ notation:
\begin{itemize}
\item[(i)] We have the energy bound $\Energy(A,B) \gtrsim K^{O(1)} \mu(A)^{3/2} \mu(B)^{3/2}$.
\item[(ii)]  There exists an open subset $E \subset A \times B$ of measure $\mu \otimes \mu(E) \gtrsim K^{O(1)} \mu(A) \mu(B)$
such that $\mu( \{ a \cdot b: (a,b) \in E \} )\lesssim  K^{O(1)} \mu(A)^{1/2} \mu(B)^{1/2}$.
\item[(iii)] There exists multiplicative sets $A'$, $B'$ contained in $A, B$ respectively such that
$\mu(A') \sim K^{O(1)} \mu(A)$, $\mu(B') \sim K^{O(1)} \mu(B)$, and $\mu(A' \cdot B') \sim K^{O(1)} \mu(A)^{1/2} \mu(B)^{1/2}$.
\item[(iv)] There exists a $O(K^{O(1)})$-approximate group $H$ of size $\mu(H) \sim K^{O(1)} \mu(A)^{1/2} \mu(B)^{1/2}$, and
elements $x, y \in G$ such that $\mu( A \cap (x \cdot H) ) \sim K^{O(1)} \mu(A)$ and 
$\mu( B \cap (H \cdot y) ) \sim K^{O(1)} \mu(B)$. 
\end{itemize}
\end{theorem}

\begin{proof}  
The implication (i) $\implies$ (iii) follows from Theorem \ref{bsg} and the trivial bounds $\mu(A') \leq \mu(A)$, $\mu(B') \leq \mu(B)$, and $\mu(A' \cdot B') \geq \mu(A')^{1/2} \mu(B')^{1/2}$.  Now we show that (iii) implies (iv).  Note that the trivial bound $\mu(A' \cdot B') \geq \max(\mu(A'), \mu(B'))$ already implies that $\mu(B) \sim K^{O(1)} \mu(A)$.  Using Theorem
\ref{energy-gleam}, we can find a $O(K^{O(1)})$-approximate group $H$ of measure $\mu(H) \sim K^{O(1)} \mu(A)$
and a finite set $X$ of cardinality $O(K^{O(1)})$ such that
$A' \subseteq X \cdot H$ and $B' \subseteq H \cdot X$.  By the pigeonhole principle we can thus find $x, y \in X$ such that
$\mu(A' \cap (x \cdot H) ) \gtrsim K^{O(1)} \mu(A)$ and 
$\mu(B' \cap (H \cdot y) ) \gtrsim K^{O(1)} \mu(B)$.  Since we have the trivial bounds $\mu(A' \cap (x \cdot H)) \leq \mu(A)$
and $\mu(B' \cap (H \cdot y)) \leq \mu(B)$, the claim (iv) follows.

Next, we show that (iv) implies (ii).  If we set $E := (A \cap (x \cdot H)) \times (B \cap (H \cdot y))$, then
we have the desired lower bound on $\mu \otimes \mu(E)$, and we have the upper bound
$$ \mu( \{ a \cdot b: (a,b) \in E \} ) \leq \mu( (x \cdot H) \cdot (H \cdot y) ) = \mu(H^2) \sim K^{O(1)} \mu(A)^{1/2} \mu(B)^{1/2}$$
as desired.

Finaly we show that (ii) implies (i).  If we let $C := \{ a \cdot b: (a,b) \in E \}$ then we have
$$ \int_C 1_A * 1_B(x)\ d\mu(x) \gtrsim K^{O(1)} \mu(A) \mu(B)$$
and (i) easily follows from the Cauchy-Schwarz inequality.
\end{proof}

\section{Metric entropy analogues}\label{entropy-sec}

In some applications involving non-discrete groups (e.g. Lie groups), it is not the measure or cardinality of 
a set which is of interest, but rather its entropy with respect to a metric.  

\begin{definition}[Metric entropy]  Let $X$ be a metric space and $\eps > 0$.  The \emph{metric entropy} (or \emph{Kolmogorov entropy}) $\N_\eps(X)$ is defined to be the least number of open balls of radius $\eps$ needed to cover $X$.
\end{definition}

\begin{remark}\label{entropy-equiv} There are several other formulations of metric entropy which are essentially equivalent to each other.  For instance, it is easy to see that the largest $\eps$-separated subset of $X$ has cardinality between $\N_{\eps}(X)$ and $\N_{\eps/2}(X)$.  Similarly, if $X$ is a subspace of a larger metric space $Y$, one can easily check that
the number of open balls of radius $\eps$ in $Y$ needed to cover $X$ lies between
$\N_{2\eps}(X)$ and $\N_\eps(X)$.  We shall shortly impose a volume doubling condition which will imply that $\N_\eps(X)$ and $\N_{2\eps}(X)$ are comparable in magnitude, and so we will not need to distinguish between these slightly different concepts of entropy.
\end{remark}

In order for the sum set theory to extend to metric entropy, we need some mild compatibility conditions between the metric structure,
the group structure, and the measure structure.  We axiomatize these as follows.

\begin{definition}[Reasonable metrics]  We say that a multiplicative group $G = (G,d_G)$ equipped with a metric $d_G$ 
a is \emph{locally reasonable metric group} if the following properties hold:
\begin{itemize}
\item[(i)] The topology on $G$ is compatible with the metric $d_G$ (thus the open balls in $d_G$ generate the topology).  Also, we assume that all closed balls are compact (thus $G$ is locally compact).
\item[(ii)] The group operations are locally Lipschitz continuous. More precisely, for every compact set $\K \subseteq G$ we have the estimates
$$ d_G(g \cdot x, g \cdot y), d_G(x \cdot g, y \cdot g), d(x^{-1},y^{-1}) \sim_{\K,G} d(x,y)
$$
for all $x,y,g \in \K$.  
\item[(iii)] We have local volume doubling.  More precisely, for any $R > 0$ we have
$$ \mu(\B(1,2r)) \sim_{R,G} \mu(\B(1,r))$$
for all $0 < r < R$, where $\B(x,r)$ is the open metric ball of radius $r$ centred at $x$.
\end{itemize}
If the implied constants in the $\sim_{\K,G}$ and $\sim_{R,G}$ notation can be chosen to be independent of $\K$ and $R$ (but still dependent on the group $G$), we say that the group is \emph{globally reasonable}.
\end{definition}

\begin{examples}  The Euclidean space $\R^d$ with the usual metric and additive group structure is globally reasonable.  Any compact Lie group
with a smooth Riemannian metric will also be globally reasonable.  If one metric is locally (resp. globally) reasonable, then any other metric bilipschitz equivalent to it will also be locally (resp. globally) reasonable.
If a locally reasonable metric group is compact, then it is automatically globally reasonable.
If $G$ is a group of linear transformations on a finite-dimensional normed vector space (with the usual topology), 
then the operator norm metric $d_G(x,y) := \|x-y\|_{\operatorname{op}}$ is locally reasonable (and thus globally reasonable, if $G$ is compact).  On the other hand, groups such as $SL_2(\R)$ will not support any globally reasonable metric, due to the non-compact nature of the conjugacy classes.
\end{examples}

As we shall shortly see, when the metric is locally reasonable, all bounded sets have finite metric entropy for each $\eps > 0$.  In this paper we shall be concerned with the bounded-dimensional regime, in which we allow all constants to depend on the implied constants in the $\sim_{\K,G}$ and $\sim_{R,G}$ notation appearing in the above definition.  The issue of precise behaviour of constants on the dimension (and on other characteristics of the group)
in the high-dimensional regime is an interesting one, but we will not pursue it here.

With the above assumptions on the metric, the metric entropy can be estimated accurately by the measure of various sets.

\begin{lemma}[Multiplicative structure of balls]\label{ball}  Let $G$ be a locally reasonable metric group.  Let 
$\K$ be a compact subset of $G$, and let $R > 0$ be a compact set.  
\begin{itemize}
\item[(i)] (Approximate normality of balls) There exists constants $0 < c_{\K,R,G} < C_{\K,R,G} < \infty$ such that we have the inclusions
$$ X \cdot \B(1, c_{\K,R,G} \eps) \subseteq \B(1,\eps) \cdot X \subseteq X \cdot \B(1,C_{\K,R,G}\eps)$$
and
$$ \B(1,c_{\K,R,G}\eps) \cdot X \subseteq X \cdot \B(1,\eps) \subseteq X \cdot \B(1,C_{\K,R,G}\eps)$$
for any $X \subseteq \K$ and $0 < \eps < R$.

\item[(ii)] (Doubling property) For any $x \in \K$, $0 < \eps < R$, and $A > 0$, we have
$$ \mu( \B(x,A\eps) ) \sim_{A,\K,R,G} \mu( \B(1, \eps) )$$

\item[(iii)] (Self-covering property) For any $x \in \K$, $0 < \eps < R$, and $A > 0$, we can cover $\B(x,A\eps)$ by
$O_{A,\K,R,G}(1)$ balls of radius $\eps$.
\end{itemize}
If $G$ is globally reasonable, then we can omit the dependence on $\K$ and $R$ in the above estimates.
\end{lemma}

\begin{proof}  If $x \in X \subseteq K$ and $0 < \eps < R$, 
then $\B(x,\eps)$ is contained in a compact set $\tilde \K = \tilde \K_{\K,R}$
which is independent of $x$ and $\eps$.  From the locally Lipschitz property we then have $d_G(x,y) \sim_{\K,R,G} d_G(1, x^{-1} \cdot y) \sim_{\K,R,G} d_G(1, y \cdot x^{-1})$ for all $x \in X$ and $y \in \B(x,\eps)$, and the claim (i) follows.

In view of (i), we see that to prove (ii) it suffices to do so when $x=1$.  But then this follows by iterating the volume doubling
property.

Finally, we prove (iii).  Let $S$ be a maximal $\eps$-separated subset of $\B(x,A\eps)$.  Clearly the balls $\B(s,\eps)$ with $s \in S$ cover $\B(x,A\eps)$.  Also, the balls $\B(s,\eps/2)$ with
$s \in S$ are disjoint subsets of $\B(x,(A+1/2)\eps)$, and thus 
$$ \sum_{s \in S} \mu( \B(s,\eps/2) ) \leq \B(x,(A+1/2)\eps).$$
Applying (ii) we obtain $|S| = O_{A,\K,R,G}(1)$, and the claim follows.  (One could also have proceeded using Lemma \ref{cover}.)

If $d_G$ is globally reasonable, then an inspection of the above arguments shows that the constants which depended on $\K$ and $R$ are now uniform in those parameters.
\end{proof}

\begin{lemma}[Relationship between entropy and measure]\label{entropy-lemma}  Let $G$ be a locally reasonable metric group.  Let 
$\K$ be a compact subset of $G$, and let $R > 0$.  Then for every $0 < \eps < R$ and every $X \subseteq \K$ we have
\begin{equation}\label{neps}
 \N_\eps(X) \sim_{\K,R,G} \frac{\mu( X \cdot \B(1,\eps) )}{\mu( \B(1,\eps) )}
\end{equation}
and
\begin{equation}\label{neps2}
 \N_\eps(X) \sim_{\K,R,G} \N_{2\eps}(X).
 \end{equation}
In particular
\begin{equation}\label{boxbox}
 \mu(X \cdot \B(1,\eps)) \sim_{\K,R,G} \mu(X \cdot \B(1,2\eps) ).
 \end{equation}
Also, the ball $\B(1,\eps)$ is approximately normal in the sense that
\begin{equation}\label{boxy}
 \mu( \B(1,\eps) \cdot X \cdot Y ) \sim_{\K,R,G} \mu( X \cdot \B(1,\eps) \cdot Y ) \sim_{\K,R,G} \mu( X \cdot Y \cdot \B(1,\eps) )
 \end{equation}
for any $X,Y \subseteq \K$.  In particular we have $\mu( X \cdot \B(1,\eps) ) \sim_{\K,R,G} \mu( \B(1,\eps) \cdot X )$.

If $G$ is globally reasonable, then we can replace $\sim_{\K,R,G}$ by $\sim_G$ in the above estimates.
\end{lemma}

\begin{remark}
Note that no measurability conditions are required on $X$ and $Y$, since sets such as $X \cdot \B(1,\eps)$ are automatically
open and precompact.  These types of inequalities are well known for Euclidean space, and the assumptions we have placed on the metric will allow us to extend the Euclidean space arguments to this more general setting without difficulty.
\end{remark}

\begin{proof}  Fix $X, \eps$.  Let $0 < c = c_{\K,R,G} \leq 1$ be a small constant to be chosen later.  Let $S$ be a maximal $c\eps$-separated
subset of $X$, then the balls $\B(s,c\eps/2)$ for $s \in S$ are disjoint, and the balls $\B(s,c\eps)$ cover $X$.  In particular
$\N_\eps(X) \leq |S|$.  By Lemma \ref{ball}(i), the balls $\B(s,c\eps/2)$ are all contained in $X \cdot \B(1,\eps)$ if $c$ is sufficiently small, and thus
$$ \sum_{s \in S} \mu( \B(s,c\eps/2) ) \leq \mu(X \cdot \B(1,\eps) ).$$
Applying Lemma \ref{ball}(ii) we have
$$ \sum_{s \in S} \mu( \B(s,c\eps/2) ) \sim_{\K,R,G} |S| \mu(\B(1,\eps)) \geq \N_\eps(X) \mu(\B(1,\eps))$$
thus obtaining the upper bound in \eqref{neps}.

Now we obtain the lower bound in \eqref{neps}.  
Let $\{ \B(s,\eps): s \in S \}$ be any covering of $X$ by $\eps$-balls with $S \subseteq X$.  
Our task is to show
that $\mu(X \cdot \B(1,\eps)) = O_{\K,R,G}( |S| \mu(\B(1,\eps)) )$.  Since $X \cdot \B(1,\eps)$ is covered by $\B(s,\eps) \cdot \B(1,\eps)$, it thus suffices by the union bound to establish that
$$ \mu( \B(s,\eps) \cdot \B(1,\eps) ) = O_{\K,R,G}( \mu(\B(1,\eps)) ).$$
But from Lemma \ref{ball}(i) we have $\B(s,\eps) \cdot \B(1,\eps) \subseteq \B(s,C\eps)$ for some $C = O_{\K,R,G}(1)$, and the claim
then follows from Lemma \ref{ball}(ii).

Now we establish \eqref{neps2}.  The bound $\N_{2\eps}(X) \leq \N_{\eps}(X)$ is trivial, so it suffices to establish
the reverse bound $\N_\eps(X) = O_{\K,R,G} \N_{2\eps}(X)$.
Let $\{ \B(s,2\eps): s \in S \}$ be a covering of $X$ by balls of radius $2\eps$ for some $S \subseteq X$; we may take $|S| = \N_{2\eps}(X)$.  It will suffice to show that $X$ can be covered by $O_{\K,R,G}(|S|)$ balls of radius $\eps$ with centres in $X$.
By Remark \ref{entropy-equiv}, this will follow if we can cover $X$ by $O_{\K,R,G}(|S|)$ balls of radius $\eps/2$ whose centres do not necessarily lie in $X$.  But this follows from Lemma \ref{ball}(iii).

Finally, the claim \eqref{boxy} follows easily from Lemma \ref{ball}(i) and \eqref{boxbox}.

If $G$ is globally reasonable, then an inspection of the above arguments shows that the constants which depended on $\K$ and $R$ are now uniform in those parameters.
\end{proof}

Lemma \ref{entropy-lemma} allows us to pass back and forth between entropies and measures, after paying various normalizing factors of $\mu(\B(1,\eps))$.  Using this lemma, one can transfer\footnote{An alternate approach would be to repeat the \emph{proofs} of the previous estimates in the metric entropy setting.  That approach also works, and in fact leads to slightly better implied constants in the $O()$ notation, however the repetition of the arguments would be rather boring and we have elected instead to illustrate the transference approach.} most of the continuous estimates of preceding sections to entropy ones, though if the metric $d_G$ is merely locally reasonable instead of globally reasonable, then one has to restrict the sets in question to a fixed compact region.  We shall focus attention on the three main results of previous sections, namely
Theorems \ref{tripling-classify}, \ref{energy-gleam}, \ref{energy-gleam-2}.  We state these results for locally reasonable metric groups, but there is an obvious variant for globally reasonable metric groups in which the dependencies of the constants on $\K$ and $R$ are dropped.

\begin{theorem}\label{tripling-classify-metric}
Let $G$ be a locally reasonable metric group.  Let $\K$ be a compact set in $G$, let $R > 0$, let $0 < \eps < R$, let $K \geq 1$, and let $A \subseteq \K$ be non-empty.  
Then the following three statements are equivalent, in the sense that if one of them holds for one choice of implied constant in the $O()$ and $\lesssim$ notation, then the other statements hold for a different choice of implied constant in the $O()$ and $\lesssim$ notation:
\begin{itemize}
\item[(i)] We have the tripling bound $\N_\eps(A^3) \lesssim_{\K,R,G} K^{O(1)} \N_\eps(A)$.
\item[(ii)] We have $\N_\eps(A^{\epsilon_1} \ldots A^{\epsilon_n} ) \sim_{\K,R,G,n} K^{O_n(1)} \N_\eps(A)$ for all $n \geq 1$ and all signs $\epsilon_1, \ldots, \epsilon_n \in \{-1,1\}$.
\item[(iii)] There exists a $O_{\K,R,G}(K^{O(1)})$-approximate group $H$ of with $\N_\eps(H) \sim_{\K,R,G} K^{O(1)} \N_\eps(A)$ 
which contains $A$, and is contained in a compact set $\tilde \K = \tilde \K(\K,R)$ depending only on $\K$ and $R$.
\end{itemize}
\end{theorem}

\begin{proof}  Let us first prove that (iii) implies (i).  By Lemma \ref{ball}(i) we have
$$ A^3 \subseteq H^3 \subseteq X \cdot X \cdot H$$
where $X$ is the set of cardinality at most $O_{\K,R,G}(K^{O(1)})$ associated to $H$.  
Using Lemma \ref{entropy-lemma} we conclude that
\begin{align*}
\N_\eps(A^3) &\lesssim_{\K,R,G} \mu( X \cdot X \cdot H \cdot \B(1,\eps) ) / \mu(\B(1,\eps)) \\
&\lesssim_{\K,R,G} |X|^2 \mu(H \cdot \B(1,\eps)) / \mu(\B(1,\eps)) \\
&\lesssim_{\K,R,G} K^{O(1)} \N_\eps(H) \\
&\sim_{\K,R,G} K^{O(1)} \N_\eps(A)
\end{align*}
which is (i).

A similar argument shows that (iii) implies (ii) and is left to the reader.  Since (ii) trivially implies (i), it remains
to show that (i) implies (iii).  
From the hypothesis on $A$ and Lemma \ref{entropy-lemma} we have
$$ \mu( A^3 \cdot \B(1,\eps) ) \lesssim_{\K,R,G} K^{O(1)} \mu(A \cdot \B(1,\eps) ).$$
Applying many applications of \eqref{boxy}, \eqref{boxbox} we conclude that
$$ \mu( (A \cdot \B(1,\eps))^3 ) \lesssim_{\K,R,G} K^{O(1)} \mu(A \cdot \B(1,c\eps) ).$$
By Theorem \ref{tripling-classify} there exists a $O_{\K,R,G}(K^{O(1)})$-approximate group $H$
which contains $A \cdot \B(1,c\eps)$, and which obeys the estimate
$$ \mu(H) \lesssim_{\K,R,G} K^{O(1)} \mu(A \cdot \B(1,\eps)) \sim_{\K,R,G} K^{O(1)} \N_\eps(A) \mu(\B(1,\eps)).$$
From the proof of Theorem \ref{tripling-classify}, and the hypothesis that $A \subseteq \K$ and $0 < \eps < R$, we also see that
$H \subseteq \tilde K$ for some compact $\tilde \K = \tilde \K(\K,R)$.  Then by Lemma \ref{entropy-lemma}
\begin{align*}
\N_\eps(H) &\lesssim_{\K,R,G} \mu( \B(1,\eps) \cdot H ) / \mu(\B(1,\eps))\\
&\leq \mu(H \cdot H) / \mu(\B(1,\eps)) \\
&\lesssim_{\K,R,G} K^{O(1)} \mu(H) / / \mu(\B(1,\eps)) \\
&\lesssim_{\K,R,G} K^{O(1)} \N_\eps(A)
\end{align*}
and (iii) follows.
\end{proof}

Now we give the metric entropy analogue of Theorem \ref{energy-gleam}.

\begin{theorem}\label{energy-gleam-metric}  
Let $G$ be a locally reasonable metric group.  Let $\K$ be a compact set in $G$, let $0 < \eps < R$, let $K \geq 1$, and let $A,B \subseteq \K$ be non-empty.  
Then the following two statements are equivalent, in the sense that if one of them holds for one choice of implied constant in the $O()$ and $\lesssim$ notation, then the other statement holds for a different choice of implied constant in the $O()$ and $\lesssim$ notation:
\begin{itemize}
\item[(i)] We have the product bound $\N_\eps(A \cdot B) \lesssim_{\K,R,G} K^{O(1)} \N_\eps(A)^{1/2} \N_\eps(B)^{1/2}$.
\item[(ii)] There exists a $O_{\K,R,G}(K^{O(1)})$-approximate group $H$ with $\N_\eps(H) \lesssim_{\K,R,G} K^{O(1)} \N_\eps(A)^{1/2} \N_\eps(B)^{1/2}$ and a finite set $X$ of cardinality $O_{\K,R,G}(K^{O(1)})$ such that $A \subset X \cdot H$ and $B \subset H \cdot X$.  Furthermore, $H$ and $X$ lie in a compact set $\tilde K = \tilde \K(\K,R)$ depending only on $\K$ and $R$.
\end{itemize}
\end{theorem}

\begin{proof}  First we show that (ii) implies (i).  We have a set $Y \subset H^2$ of cardinality $O_{\K,R,G}(K^{O(1)})$
such that $H \cdot H \subseteq H \cdot Y$, which implies that $A \cdot B \subseteq X \cdot H \cdot Y \cdot X$.
From Lemma \ref{entropy-lemma} we then have
\begin{align*}
\N_\eps( A \cdot B ) &\leq \N_\eps( X \cdot H \cdot Y \cdot X )\\
&\lesssim_{\K,R,G} \mu( X \cdot H \cdot Y \cdot X \cdot \B(1,\eps) ) / \mu(\B(1,\eps)) \\
&\lesssim_{\K,R,G} |X|^2 |Y| \mu( H \cdot \B(1,\eps) ) / \mu(\B(1,\eps)) \\
&\lesssim_{\K,R,G} K^{O(1)} \N_\eps(H) \\
&\lesssim_{\K,R,G} K^{O(1)} \N_\eps(A)^{1/2} \N_\eps(B)^{1/2}
\end{align*}
which is (i).

Now we show that (i) implies (ii).   From the hypothesis and Lemma \ref{entropy-lemma} we have
\begin{align*}
\mu( A \cdot \B(1,\eps) \cdot B \cdot \B(1,\eps) ) &\lesssim_{\K,R,G} \mu(A \cdot B \cdot \B(1,\eps) ) \\
&\lesssim_{\K,R,G} \N_\eps(A \cdot B) \mu(\B(1,\eps))\\
&\lesssim_{\K,R,G} K^{O(1)} \N_\eps(A)^{1/2} \N_\eps(B)^{1/2} \mu(\B(1,\eps)) \\
&\sim_{\K,R,G} K^{O(1)} \mu(A \cdot \B(1,\eps))^{1/2} \mu(B \cdot \B(1,\eps))^{1/2}.
\end{align*}
Applying Theorem \ref{energy-gleam}, we can find a $O_{\K,R,G}(K^{O(1)})$-approximate group $H$ and a set $X$ such that
$A \cdot \B(1,\eps) \subset X \cdot H$ and $B \cdot \B(1,\eps) \subset H \cdot X$, with the bounds
\begin{align*}
\mu(H) &\lesssim_{\K,R,G} K^{O(1)} \mu(A \cdot \B(1,\eps))^{1/2} \mu(B \cdot \B(1,\eps))^{1/2} \\
&\lesssim_{\K,R,G} \N_\eps(A)^{1/2} \N_\eps(B)^{1/2} \mu(\B(1,\eps)).
\end{align*}
and $|X| \lesssim_{\K,R,G} K^{O(1)}$.  We then have
\begin{align*}
\N_\eps(H) &\lesssim_{\K,R,G} \mu(\B(1,\eps) \cdot H) / \mu(\B(1,\eps)) \\
&\leq \mu( A \cdot \B(1,\eps) \cdot H ) / \mu(\B(1,\eps)) \\
&\leq \mu( X \cdot H \cdot H ) / \mu(\B(1,\eps)) \\
&\lesssim_{\K,R,G} |X| K^{O(1)} \mu(H) / \mu(\B(1,\eps)) \\
&\lesssim_{\K,R,G} K^{O(1)} \N_\eps(A)^{1/2} \N_\eps(B)^{1/2}
\end{align*}
and (ii) follows.
\end{proof}

Now we turn to developing a metric entropy analogue of Theorem \ref{energy-gleam-2}.  This will be a bit trickier as we
shall need an ``$\eps$-approximate'' version of the multiplicative energy $\Energy(A,B)$.  There are a number of essentially equivalent 
ways to do so, each of which are at least somewhat artificial; for sake of concreteness we shall fix one such as follows.  Given any $A, B \subset G$ and $\eps > 0$, the set
$$ Q_\eps(A,B) := \{ (a,b,a',b') \in A \times B \times A \times B: d( a \cdot b, a' \cdot b' ) \leq \eps \}$$
of approximately multiplicative quadruples is a subset of $G^4$, which we view as a metric space with the metric
$$ d_{G^4}((x_1,x_2,x_3,x_4),(y_1,y_2,y_3,y_4)) := \sum_{i=1}^4 d_G(x_i,y_i).$$
We then define the $\eps$-approximate multiplicative energy $\Energy_\eps(A,B)$ to be the quantity $\N_\eps(Q_\eps(A,B))$.  Note
that if $A,B$ are finite sets, then this quantity will equal the usual (discrete) multiplicative energy
\eqref{energy-discrete} for $\eps$ sufficiently small.

\begin{theorem}\label{energy-gleam-2-metric}  
Let $G$ be a locally reasonable metric group.  Let $\K$ be a compact set in $G$, let $0 < \eps < R$, let $K \geq 1$, and let $A,B \subseteq \K$ be non-empty.  
Then the following four statements are equivalent, in the sense that if one of them holds for one choice of implied constant in the $O()$ and $\lesssim$ notation, then the other statement holds for a different choice of implied constant in the $O()$ and $\lesssim$ notation:
\begin{itemize}
\item[(i)] We have the energy bound $\Energy_\eps(A,B) \gtrsim_{\K,R,G} K^{O(1)} \N_\eps(A)^{3/2} \N_\eps(B)^{3/2}$.
\item[(ii)]  There exists a subset $E \subset A \times B$ of entropy $\N_\eps(E) \gtrsim_{\K,R,G} K^{O(1)} \N_\eps(A) \N_\eps(B)$
such that $\N_\eps( \{ a \cdot b: (a,b) \in E \} ) \lesssim_{\K,R,G} K^{O(1)} \N_\eps(A)^{1/2} \N_\eps(B)^{1/2}$.  (Of course, we measure the entropy of $E$ using the product metric on $G^2$.)
\item[(iii)] There exists subsets $A'$, $B'$ of $A, B$ respectively such that
$\N_\eps(A') \sim_{\K,R,G} K^{O(1)} \N_\eps(A)$, $\N_\eps(B') \sim_{\K,R,G} K^{O(1)} \N_\eps(B)$, and $\N_\eps(A' \cdot B') \sim_{\K,R,G} K^{O(1)} \N_\eps(A)^{1/2} \N_\eps(B)^{1/2}$.
\item[(iv)] There exists a $O(K^{O(1)})$-approximate group $H$ with $\N_\eps(H) \sim_{\K,R,G} K^{O(1)} \N_\eps(A)^{1/2} \N_\eps(B)^{1/2}$, and
elements $x, y \in G$ such that $\N_\eps( A \cap (x \cdot H) ) \sim_{\K,R,G} K^{O(1)} \N_\eps(A)$ and 
$\N_\eps( B \cap (H \cdot y) ) \sim_{\K,R,G} K^{O(1)} \N_\eps(B)$.  Furthermore, $H$, $x$, $y$ lie in a compact set $\tilde \K = \K(\K,R)$ that depends only on $\K$ and $R$.
\end{itemize}
\end{theorem}

\begin{proof}  Let us first show that (iv) implies (iii).  We set $A' := A \cap (x \cdot H)$ and $B' := B \cap (H \cdot y)$.
The bounds on $\N_\eps(A')$ and $\N_\eps(B')$ are obvious.  From Lemma \ref{entropy-lemma} and the trivial estimate
$\mu( X \cdot Y ) \geq \mu(X)^{1/2} \mu(Y)^{1/2}$ we see that $\N_\eps(A' \cdot B') \gtrsim_{\K,R,G} \N_\eps(A')^{1/2} \N_\eps(B')^{1/2}$, which gives the lower bound on $\N_\eps(A' \cdot B')$.  To obtain the upper bound, we use Lemma \ref{entropy-lemma} to compute
\begin{align*}
\N_\eps(A' \cdot B') &\leq \N_\eps( x \cdot H \cdot H \cdot y ) \\
&\lesssim_{\K,R,G} \mu( x \cdot H \cdot H \cdot y \cdot \B(1,\eps) ) / \mu(\B(1,\eps)) \\
&\lesssim_{\K,R,G} \mu( \B(1,\eps) \cdot H \cdot H ) / \mu(\B(1,\eps)).
\end{align*}
But $H \cdot H$ is covered by $O_{\K,R,G}(K^{O(1)})$ right-translates of $H$, and so
$$ \N_\eps(A' \cdot B')  \lesssim_{\K,R,G} K^{O(1)} \N_\eps(H) \sim_{\K,R,G} K^{O(1)} \N_\eps(A)^{1/2} \N_\eps(B)^{1/2}$$
as desired.

Now we show that (iii) implies (ii).  We take $E := A' \times B'$.  The bound on $\N_\eps( \{ a \cdot b: (a,b) \in E \} )$
is obvious, while by considering products of $\eps$-separated sets it is easy to establish a bound of the form
$$\N_{\eps/100}(E) \gtrsim_{\K,R,G} \N_\eps(A') \N_\eps(B') \gtrsim_{\K,R,G} K^{O(1)} \N_\eps(A) \N_\eps(B).$$
The claim then follows from \eqref{neps2} (note that if $G$ is locally reasonable then so is $G^2$).

Now we show that (ii) implies (i).  Let $E'$ be a maximal $100\eps$-separated subset of $E$, then by \eqref{neps2}
$$ |E'| \geq \N_{100\eps}(E) \sim_{\K,R,G} \N_\eps(E) \gtrsim_{\K,R,G} K^{O(1)} \N_\eps(A) \N_\eps(B).$$
Let $D$ be a maximal $\eps/2$-separated subset of $\{ a \cdot b: (a,b) \in E \}$, thus
$$ |D| \leq \N_{\eps/4}( \{ a \cdot b: (a,b) \in E \} ) 
\sim_{\K,R,G} \N_\eps( \{ a \cdot b: (a,b) \in E \} ) \lesssim_{\K,R,G} K^{O(1)} \N_\eps(A)^{1/2} \N_\eps(B)^{1/2}.$$
Observe that for every $(a,b) \in E'$, the product $a \cdot b$ lies within $c\eps$ of an element of $D$, thus
$$ \sum_{x \in D} | \{ (a,b) \in E': d_G( a \cdot b, x ) \leq \eps/2 \} | \geq |E'| \gtrsim_{\K,R,G} K^{O(1)} \N_\eps(A) \N_\eps(B).$$
Applying Cauchy-Schwarz we conclude that
$$ \sum_{x \in D} | \{ (a,b) \in E': d_G( a \cdot b, x ) \leq \eps/2 \} |^2 \gtrsim_{\K,R,G} K^{O(1)} \N_\eps(A)^{3/2} \N_\eps(B)^{3/2}.$$
Observe that if $(a,b), (a',b') \in E'$ and $x \in D$ are such that $d_G(a \cdot b,x),d_G(a' \cdot b',x) \leq \eps/2$
then $(a,b,a',b') \in Q_\eps(A,B)$.  Thus
$$ |Q_\eps(A,B) \cap (E' \times E')| \gtrsim_{\K,R,G} K^{O(1)} \N_\eps(A)^{3/2} \N_\eps(B)^{3/2}.$$
But $E' \times E'$ is clearly $\eps$-separated, thus
$$ \Energy_\eps(A,B) = \N_\eps( Q_\eps(A,B) ) \gtrsim_{\K,R,G} K^{O(1)} \N_\eps(A)^{3/2} \N_\eps(B)^{3/2}$$
as desired.

Finally, we show that (i) implies (iv), which is the most difficult implication.  
Let $C = C_{\K,R,G}$ be a large constant to be chosen later.  Let
$\overline{A} := A \cdot \B(1,C\eps)$ and $\overline{B} := B \cdot \B(1,C\eps)$, thus from
Lemma \ref{entropy-lemma} we see that $\overline{A}$, $\overline{B}$ are multiplicative sets with
\begin{equation}\label{overbite}
\mu(\overline{A}) \sim_{\K,R,G,C} \N_\eps(A) \mu(\B(1,\eps)); \quad
\mu(\overline{B}) \sim_{\K,R,G,C} \N_\eps(B) \mu(\B(1,\eps)).
\end{equation}

Now consider the quantity $\Energy( \overline{A}, \overline{B} )$.  We can rewrite this as
\begin{align*}
\Energy( \overline{A}, \overline{B} ) &= \int_G (1_{\overline{A}} * 1_{\overline{B}})(x)^2\ d\mu(x) \\
&= \int_{\overline{A}} \int_{\overline{B}} 1_{\overline{A}} * 1_{\overline{B}}( a \cdot b)\ d\mu(a) d\mu(b) \\
&= \int_{\overline{A}} \int_{\overline{B}} \int_{\overline{A}} 1_{\overline{B}}( (a')^{-1} \cdot a \cdot b)\ d\mu(a') d\mu(a) d\mu(b).
\end{align*}
Now observe that if $C$ is large enough, we see that for any $x \in G$, the set
$B \cdot \B(1,\sqrt{C}\eps)$ intersects $\B(x,\sqrt{C}\eps)$ only when $x \in \overline{B}$.  This (and Lemma \ref{ball}(ii))
leads to the pointwise estimate
$$ 1_{\overline{B}}( x ) \gtrsim_{\K,R,G} \frac{1}{\mu(\B(1,\eps))} \int_{B \cdot \B(1,\sqrt{C}\eps)}\ 1_{\B(x,\sqrt{C}\eps)}(b) d\mu(b)$$
and hence
$$ \Energy( \overline{A}, \overline{B} ) \gtrsim_{\K,R,G} \frac{1}{\mu(\B(1,\eps))} 
\int_{\overline{A}} \int_{\overline{B}} \int_{\overline{A}} \int_{B \cdot \B(1,\sqrt{C}\eps)}
1_{\B((a')^{-1} \cdot a \cdot b,\sqrt{C}\eps)}(b')\ d\mu(b') d\mu(a') d\mu(a) d\mu(b).$$
Now observe (from the local Lipschitz property) that if $(a,b,a',b') \in Q_\eps(A,B) \cdot B_{G^4}(1,\eps)$, where $B_{G^4}(1,\eps)$
denotes the ball in $G^4$, then $d( a \cdot b, a' \cdot b' ) \lesssim_{\K,R,G} \eps$, and hence (if $C$ is large enough)
$b' \in (a')^{-1} \cdot a \cdot b,\sqrt{C}\eps)$.  Thus we have
$$ \Energy( \overline{A}, \overline{B} ) \gtrsim_{\K,R,G} \frac{1}{\mu(\B(1,\eps))} \mu^{\otimes 4}( Q_\eps(A,B) \cdot B_{G^4}(1,\eps) )$$
and hence (by the analogue\footnote{Here we need the easily verified fact that the direct product of finitely many locally reasonable metric groups is still locally reasonable.} of Lemma \ref{entropy-lemma} for $G^4$)
$$ \Energy( \overline{A}, \overline{B} ) \gtrsim_{\K,R,G} \frac{\mu^{\otimes 4}(B_{G^4}(1,\eps))}{\mu(\B(1,\eps))} 
\N_\eps( Q_\eps(A,B) ) \gtrsim_{\K,R,G} \mu(\B(1,\eps))^3 \Energy_\eps(A,B).$$
We henceforth fix $C$ to be a suitably large quantity depending on $\K,R,G$, and thus can omit the dependence of $C$ in
the estimates which follow.
By hypothesis on $\Energy_\eps(A,B)$ and \eqref{overbite}, we thus have
$$\Energy( \overline{A}, \overline{B} )  \gtrsim_{\K,R,G} K^{O(1)} \mu(\overline{A})^{3/2} \mu(\overline{B})^{3/2}.$$
Applying Proposition \ref{energy-gleam-2} we can thus locate a $O_{\K,R,G}(K^{O(1)})$-approximate 
group $H$ of size 
$$\mu(H) \sim_{\K,R,G} K^{O(1)} \mu(\overline{A})^{1/2} \mu(\overline{B})^{1/2}$$
and elements $x, y \in G$ such that $\mu( \overline{A} \cap (x \cdot H) ) \sim_{\K,R,G} K^{O(1)} \mu(\overline{A})$ and 
$\mu( \overline{B} \cap (H \cdot y) ) \sim_{\K,R,G} K^{O(1)} \mu(\overline{B})$.   An inspection of the proof of Proposition \ref{energy-gleam-2} reveals that $H$, $x$, $y$ are also contained in a compact set that depends only on $\K$ and $R$.
Note that from the trivial bounds $\mu( \overline{A} \cap (x \cdot H)) \leq \mu(\overline{A})$ and
$\mu( \overline{B} \cap (H \cdot y) ) \leq \mu(\overline{B})$ we can conclude that $\overline{A}$ and $\overline{B}$ are comparable in size:
$$ \mu( \overline{A} ) \sim_{\K,R,G} K^{O(1)} \mu(\overline{B}).$$
From \eqref{overbite} we thus have entropy comparability also:
\begin{equation}\label{entropy-compare}
\N_\eps(A) \sim_{\K,R,G} K^{O(1)} \N_\eps(B).
\end{equation}
From \eqref{overbite} again, we have a good bound on the measure of $H$:
\begin{equation}\label{muh}
\mu(H) \sim_{\K,R,G} K^{O(1)} \N_\eps(A)^{1/2} \N_\eps(B)^{1/2} \mu(\B(1,\eps)).
\end{equation}
However to get a good bound on the \emph{entropy} of $H$ we need to estimate $\mu(H \cdot B(1,\eps))$.  This we shall do by
means of the Ruzsa triangle inequality.  Observe that
$$ \mu( [\overline{A} \cap (x \cdot H)] \cdot H ) \leq \mu( x \cdot H \cdot H ) \lesssim_{\K,R,G} K^{O(1)} \mu(H)$$
since $H$ is an approximate group.  From our bounds on $\mu(H)$ and $\mu(\overline{A} \cap (x \cdot H))$
we conclude that
$$ d( \overline{A} \cap (x \cdot H), H^{-1} ) \leq O(\log K) + O_{\K,R,G}(1).$$
Next, we observe that
$$ \mu( [\overline{A} \cap (x \cdot H)] \cdot \B(1,\eps) ) \leq \mu( \overline{A} \cdot \B(1,\eps) )
= \mu( A \cdot \B(1,C\eps) \cdot \B(1,\eps) )$$
and so from Lemma \ref{entropy-lemma} we have
$$ \mu( [\overline{A} \cap (x \cdot H)] \cdot \B(1,\eps) ) \lesssim_{\K,R,G} K^{O(1)} \N_\eps(A) \mu(\B(1,\eps)).$$
This gives a bound on the Ruzsa distance:
$$ d( \overline{A} \cap (x \cdot H), \B(1,\eps)^{-1} ) \leq \log \N_\eps(A) + O(\log K) + O_{\K,R,G}(1).$$
By the triangle inequality we conclude that
$$ d( H^{-1}, \B(1,\eps)^{-1} ) \leq \log \N_\eps(A) + O(\log K) + O_{\K,R,G}(1)$$
which implies that
\begin{equation}\label{buh}
\mu(H \cdot \B(1,\eps)) \lesssim_{\K,R,G} K^{O(1)} \N_\eps(A) \mu(\B(1,\eps)).
\end{equation}
Combining this with Lemma \ref{entropy-lemma}, \eqref{muh}, \eqref{entropy-compare},
and the trivial lower bound $\mu(H \cdot \B(1,\eps)) \geq \mu(H)$ we conclude
that
$$ \N_\eps(H) \sim_{\K,R,G} K^{O(1)} \N_\eps(A)^{1/2} \N_\eps(B)^{1/2}.$$
We are nearly done with establishing (iv), but there is one slight problem: we have shown that $\overline{A}$ and $\overline{B}$ have large intersection (in the measure sense) with translates of $H$, but we need $A$ and $B$ to have large intersection (in the entropy sense) with translates of $H$.  There are a number of ways to resolve this; one is as follows.  Observe that
$$ \N_\eps( \overline{A} \cap (x \cdot H) ) \gtrsim_{\K,R,G} \mu( \overline{A} \cap (x \cdot H ) ) / \mu(\B(1,\eps))
\gtrsim_{\K,R,G} K^{O(1)} \N_\eps(A).$$
Thus there exists a $\eps$-separated subset $\overline{A'}$ of $\overline{A} \cap (x \cdot H)$ of cardinality
$\gtrsim_{\K,R,G} K^{O(1)} \N_\eps(A)$.  Using \eqref{neps2} one can refine this subset to be $C\eps$-separated for any fixed $C = O_{\K,R,G}(1)$ without degrading the cardinality of $\overline{A'}$ significantly.  By construction of $\overline{A'}$, each element of $\overline{A'}$ is at a distance $O_{\K,R,G}(\eps)$ to an element of $A$.  This shows that there exists an $\eps$-separated
subset $A'$ of $A$ of cardinality $\gtrsim_{\K,R,G} K^{O(1)} \N_\eps(A)$, with each element of $A'$ at a distance $O_{\K,R,G}(\eps)$ to an element of $x \cdot H$.  If we then set $\tilde H := H \cdot \B(1,C\eps)$ for some sufficiently large $C = O_{\K,R,G}(1)$,
we see that $A'$ is contained in $x \cdot \tilde H$ and so $\N_\eps(A \cap (x \cdot \tilde H)) \gtrsim_{\K,R,G} K^{O(1)} \N_\eps(A)$.
A similar argument yields $\N_\eps((y \cdot \tilde H) \cap B) \gtrsim_{\K,R,G} K^{O(1)} \N_\eps(B)$.  Now we need to pass from $\tilde H$ back to $H$.  Recall that as $H$ is an approximate group, $H \cdot H$ can be covered by $O_{\K,R,G}(K^{O(1)})$ left-translates of $H$, hence $H \cdot \tilde H$ can be covered by $O_{\K,R,G}(K^{O(1)})$ left-translates of $\tilde H$.
Combining this with \eqref{muh}, \eqref{buh}, \eqref{boxbox} we have
$$ \mu(H \cdot \tilde H) \lesssim_{\K,R,G} K^{O(1)} \mu(H).$$
Applying Lemma \ref{cover} we can thus cover $\tilde H$ by $O_{\K,R,G}(K^{O(1)})$ left-translates of $H$.  In particular
we can cover $x \cdot \tilde H$ by $O_{\K,R,G}(K^{O(1)})$ sets of the form $x' \cdot H$, and so by the pigeonhole principle
we have $\N_\eps(A \cap (x' \cdot H)) \gtrsim_{\K,R,G} K^{O(1)} \N_\eps(A)$ for some $x'$, which one can easily verify is contained in a compact set $\tilde \K(\K,R)$ depending only on $\K$ and $R$.  A similar argument (using \eqref{boxy} to move the $\B(1,\eps)$ factors around as necessary) gives $\N_\eps(B \cap (H \cdot y')) \gtrsim_{\K,R,G} K^{O(1)} \N_\eps(B)$ for some $y' \in \tilde \K(\K,R)$, and (iv) follows.
\end{proof}

One can of course develop metric entropy analogues of many of the other estimates from previous sections (such as the Ruzsa triangle inequality).  We leave the details to the reader.

\section{Inverse theorems}\label{inverse-sec}

The above theory reduces the study of sets of small doubling or tripling (or pairs of sets with small product set, small partial product set, or large multiplicative energy) to that of studying approximate groups, at least if one is prepared to lose polynomial factors in the constants and (in the locally reasonable metric entropy setting) one restricts all sets to a compact region.  There remains of course the question of how to effectively classify these approximate groups; we refer to this as the \emph{inverse product set problem} (or the \emph{inverse sum set problem}, in the abelian additive setting).  At present, there is not even a reasonable conjecture as to what such objects should look like; there are obvious examples of approximate groups, such as genuine groups, geometric progressions\footnote{In the abelian case, the group $G$ is usually written additively, and it is then the \emph{arithmetic} progressions which are relevant here.  However as we are considering the non-commutative setting we are forced to depart from the usual additive notation and work instead with geometric progressions.}, and (given sufficient commutativity) the direct sum of other approximate groups, but it is not clear in general what the statement should be\footnote{Indeed, the problem can be viewed as a robust version of the problem of classifying all the subgroups of a given group $G$, which is already quite a difficult problem, especially for highly non-abelian groups such as the permutation group $S_n$.  In some cases it seems that the class of approximate subgroups of $G$ is not much ``richer'' the class of genuine subgroups of $G$, in the sense that one can express approximate subgroups as dense subsets of combinations of genuine subgroups and related objects such as geometric progressions, but this might not be true for sufficiently complicated groups $G$.}.

There are however a number of special cases which are well understood.  If $G$ is a discrete abelian $r$-torsion group for some small
$r > 1$ (thus $x^r = 1$ for all $r$) and $A$ is a finite non-empty subset of $G$ 
then it is known that $|A \cdot A| = O(|A|)$ if and only if $A$ can be contained in a finite subgroup $H$ of $G$ of size $O(|A|)$; see \cite{ruzsa-group}.  If $G$ is instead a discrete abelian torsion-free group, and $A$ is a finite non-empty subset of $G$, then a famous theorem of Freiman \cite{frei} (see also \cite{bilu-freiman}, \cite{ruzsa-freiman}, \cite{chang}) shows that $|A \cdot A| = O(|A|)$ if and only if $A$ is contained in the product $P$ of $O(1)$ geometric progressions, whose total cardinality is $O(|A|)$.  These results were unified in \cite{gr-4}, in which $A$ was now a finite non-empty subset of an arbitrary abelian group $G$, and the result now being that $|A \cdot A| = O(|A|)$ if and only if $A$ is contained in the product $P$ of $O(1)$ geometric progressions and a finite subgroup of $G$, whose total cardinality is $O(|A|)$.  See \cite{tao-vu} for a presentation of all of these abelian results.  Apart from the (important) issue of quantifying the dependence of constants here, this is a satisfactory resolution of the inverse sum set problem in the discrete setting.

In the abelian setting it is also easy to pass to the continuous setting and the metric entropy setting.  For instance, we have

\begin{proposition}[Continuous version of Freiman's theorem]\label{frei-cts} Let $d \geq 1$.  Let $A$ be an open bounded non-empty subset of $\R^d$ such that
$\mu(A + A) \leq K \mu(A)$ for some $K \geq 2^d$, where $\mu$ denotes Lebesgue measure.  Then there exists an $\eps > 0$ and a set $P$ which is the sum of $O_K(1)$ arithmetic progressions in $\R^d$ such that $A \subseteq P + \B(0,\eps)$ and $\mu(P + \B(0,\eps)) \sim_K \mu(A)$.
\end{proposition}

\begin{remark}\label{cfr} Note that the trivial inclusion $A+A \supset 2 \cdot A$ (or the Brunn-Minkowski inequality) shows that $K$ cannot be less than $2^d$.  In the converse direction, it is easy to see that if $A \subseteq P + \B(0,\eps)$ and $\mu(P + \B(0,\eps)) \sim_K \mu(A)$, where $P$ is the sum of $O(1)$ arithmetic progressions, then $\mu(A+A) \sim_K \mu(A)$.  One can certainly
use the arguments in \cite{bilu-freiman}, \cite{ruzsa-freiman}, \cite{gt-freimanbilu} to quantify the exact dependence on $K$ in the above proposition but we will not attempt to do so here.  It is also not difficult to modify the above proposition to replace the Euclidean space $\R^d$ with a torus such as $\R^d/\Z^d$ by a lifting argument; we omit the details.
\end{remark}

\begin{proof}  Since $A+A$ is open, we have $A+A = \bigcap_{\eps > 0} A+A + \B(0,\eps)$.  By the monotone convergence theorem, we can
thus find an $\eps > 0$ such that $\mu( A + A + B(0,\eps) ) \sim \mu( A + A ) \sim_K \mu( A )$.  

Now let $\tilde A := (A + B(0,\eps/2)) \cap ( \frac{\eps}{10d} \cdot \Z^d )$, thus $A$ is a finite non-empty set.  
From the inclusions
$$ A \subseteq \tilde A + B(0,\eps) \hbox{ and } \tilde A + \tilde A + B(0, \frac{\eps}{100d} ) \subseteq A + A + B(0,\eps)$$
one easily verifies the estimate
$$ \mu(A) \lesssim_K |\tilde A| \eps^d \leq |\tilde A + \tilde A| \eps^d \lesssim_{K} \mu(A).$$
Note that any dependencies on $d$ of the implied constant can be converted to a dependency on $K$ since $K \geq 2^d$.  In particular
we have $|\tilde A + \tilde A| \sim_K |\tilde A|$.  Applying Freiman's theorem (see e.g. \cite{frei}, \cite{bilu-freiman}, \cite{ruzsa-freiman}, \cite{chang},  \cite{nathanson}, \cite{tao-vu}) we can thus place $\tilde A$ inside a set $P \subset \R^d$ of cardinality
$|P| \sim_K |\tilde A|$ which is the sum of $O_K(1)$ arithmetic progressions.  Since $A \subseteq \tilde A + B(0,\eps) \subseteq P + B(0,\eps)$, we have
$$ \mu(A) \leq \mu( P + B(0,\eps) ) \lesssim_d |P| \eps^d \sim_K |\tilde A| \eps^d \sim_K \mu(A)$$
and the claim follows.
\end{proof}

\begin{proposition}[Entropy version of Freiman's theorem]\label{frei-entropy} Let $d \geq 1$ and $\eps > 0$.  Let $A$ be a bounded non-empty subset of $\R^d$ such that $\N_\eps(A + A) \leq K \N_\eps(A)$ for some $K \geq 1$.  Then there exists a set $P$ which is the sum of $O_{K,d}(1)$ arithmetic progressions in $\R^d$ such that $A \subseteq P + B(0,\eps)$ and $|P| \sim_{K,d} \N_\eps(A)$.
\end{proposition}

\begin{remark} The commentary in Remark \ref{cfr} also applies in this setting.  For instance, if $A \subseteq P + B(0,\eps)$ and $|P| \sim_{K,d} \N_\eps(A)$ and $P$ is the sum of $O_{K,d}(1)$ arithmetic progressions then it is easy to see that $\N_\eps(A+A) \sim_{K,d} \N_\eps(A)$.  
\end{remark}

\begin{proof}  Again, we take $\tilde A := (A + B(0,\eps/2)) \cap (\frac{\eps}{10d} \cdot \Z^d)$.  From
Lemma \ref{entropy-lemma} (and the global reasonableness of $\R^d$) we 
have $|\tilde A| \sim_d \N_\eps(A)$ and $|\tilde A + \tilde A| \lesssim_d \N_\eps(A+A)$, and thus
$|\tilde A + \tilde A| \sim_{K,d} |\tilde A|$.  We now argue as in the proof of Proposition \ref{frei-cts}.
\end{proof}

We now turn to the noncommutative setting.  Here our understanding is only satisfactory for a few special noncommutative groups; in the general case it is not even clear what the correct statement of an inverse product setting theorem should be, let alone how to prove it.
We shall restrict our attention to the discrete setting, in other words in understanding
those finite non-empty sets $A$ for which $|A \cdot A| = O(|A|)$ (or $|A \cdot A \cdot A| = O(|A|)$), or for classifying finite non-empty
$O(1)$-approximate groups; in view of the preceding results it seems likely that the transferral of the discrete results to a continuous or metric entropy setting will not be too difficult.

Inverse product set theorems for groups of affine or projective mappings on the real or complex line or projective line, and hence to groups
such as $SL_2(\R)$ or $SL_2(\C)$, were studied in \cite{el1,el2,el3,elk,elru}.  A typical result here is that if $A \subset SL_2(\C)$
is a finite non-empty set such that $|A \cdot A^{-1}| = O(|A|)$, then $A$ is contained inside $O(1)$ left-cosets
of an abelian subgroup of $SL_2(\C)$.

The case $G = SL_2(\Z/p\Z)$, with $p$ a large prime, was studied in \cite{Helf}.  In particular it was shown that if
$A \subset G$ had size $p^\eps \leq |A| \leq p^{-\eps} |G|$ for some $\eps > 0$, and $A$ was not contained in any proper subgroup of $G$, then one had the tripling estimate $|A \cdot A \cdot A| \geq p^{\delta} |A|$ for some $\delta = \delta(\eps) > 0$.  Thus the only sets of small tripling are those sets which are very small, very large, or are contained in a proper subgroup (e.g. a geometric progression containing the identity).  Using the machinery in this section one can also obtain a classification of sets of small doubling, which we leave as an exercise to the reader.

Another interesting example arises in the work of Lindenstrauss \cite{linden}, in which $G$ is now the \emph{lamplighter group} $\Z \times (\Z/2\Z)^\Z$ with group law $(i,a) \cdot (j,b) = (i+j, \sigma^j a + b)$, where $\sigma$ is the standard shift on $(\Z/2\Z)^\Z$.
There it was shown that the group $G$ contains no F{\o}lner sequence of sets of small doubling constant, despite $G$ being amenable (and solvable).

The case of very small doubling, e.g. $|A \cdot A| \leq 2|A|$, was treated in \cite{kemperman}, \cite{bf}, \cite{hls} in the torsion-free non-commutative case.  In this case one has $|A \cdot A| \geq 2|A|-1$, with equality only holding when $A$ is a geometric progression.

We were unable to say anything new about the inverse product set problem for general groups.  However for discrete groups $G$ which have a normal subgroup $H$, it turns out that one can exploit the short exact sequence
$$ \{1\} \to H \to G \to G/H \to \{1\}$$
to split the inverse product set problem for $G$ into the inverse product set problem for $H$ and $G/H$ separately, together with the problem\footnote{This problem seems to be somewhat difficult, however it does appear to be fractionally simpler than the original inverse product set problem on $G$, so the reduction is not entirely trivial.  It is somewhat analogous to the reduction of the (open) ``polynomial Freiman-Ruzsa conjecture'' to a conjecture concerning approximate homomorphisms in \cite{green-survey}.} of classifying a certain type of ``approximate group homomorphism'' from an approximate subgroup of $G/H$ into $G$.  To motivate matters, let us first see how \emph{genuine} subgroups of $G$ (as opposed to approximate groups) split under this short exact sequence.

\begin{lemma}[Splitting lemma, group case]  Let $H$ be a normal subgroup of a group $G$, and let $A \subset G$.  Let $\pi: G \to G/H$ be the canonical projection.  Then the following are equivalent.
\begin{itemize}
\item[(i)] $A$ is a subgroup of $G$.
\item[(ii)] There exists a subgroup $B$ of $H$, a subgroup $C$ of $G/H$, and a partial inverse $\phi: C \to G$ to $\pi$ (i.e. $\pi(\phi(x)) = x$ for all $x \in C$) with the property
\begin{equation}\label{phib}
 \phi(x) B = B \phi(x) \hbox{ for all } x \in C
 \end{equation}
(thus $\phi$ takes values in the normaliser of $B$)
and the quotiented homomorphism property
\begin{equation}\label{phib-2}
\phi(xy) \in \phi(x) \phi(y) B \hbox{ for all } x,y \in C
\end{equation}
and such that $A$ has the representation
\begin{equation}\label{adef}
 A = \bigcup_{x \in C} \phi(x) B.
 \end{equation}
In particular, $|A| = |C| |B|$.
\end{itemize}
\end{lemma}

\begin{remark} One way to view this lemma is to think of $G$ as a principal $H$-bundle over $G/H$.  Then $A$ is a principal $B$-bundle over $C$ that takes values in the normaliser of the structure group $B$, and which collapses to a group homomorphism from $C$ to $G$ when quotiented out by $B$.
\end{remark}

\begin{proof} Let us first verify that (ii) implies (i).  It is easy to verify from \eqref{phib}, \eqref{phib-2} that $\phi(0) \in B$ and
$\phi(x^{-1}) \in \phi(x) B$ for all $x \in C$; from these facts and \eqref{phib}, \eqref{phib-2}, \eqref{adef} we quickly see 
that $A$ contains the identity, that $A^{-1} = A$, and that $A \cdot A = A$; in other words, $A$ is a subgroup of $G$.

Now let us verify that (i) implies (ii).  Let $C := \pi(A)$ and $B := A \cap K$; it is easy to see that $B$ and $C$ are subgroups of $K$ and $C/K$ respectively.  Let $\phi: C \to A$ be an arbitrary partial inverse to the map $\pi: A \to C$ (which exists thanks to the axiom of choice),
then we have \eqref{adef}.  Since $A \cdot A = A$ we conclude \eqref{phib-2}.  To verify \eqref{phib}, observe that $\phi(x) B \phi(x)^{-1}$
lies in $A \cdot A \cdot A^{-1} = A$, but also lies in the normal group $K$, and must therefore lie in $A \cap K = B$.  This shows the inclusion $\phi(x) B \subseteq B \phi(x)$, and the other inclusion is proven similarly.
\end{proof}

We now present an analogue of the above lemma for approximate groups.

\begin{lemma}[Splitting lemma, approximate group case]\label{split-approx}  Let $H$ be a normal subgroup of a discrete multiplicative group $G$, and let $\pi: G \to G/H$ be the canonical homomorphism.  Note that $H$ and $G/H$ are also discrete multiplicative groups.  Let $A \subset G$, and let $K \geq 1$.
Then the following three statements are equivalent, in the sense that if one of them holds for one choice of implied constant in the $O()$ and $\lesssim$ notation, then the other statement holds for a different choice of implied constant in the $O()$ and $\lesssim$ notation:
\begin{itemize}
\item[(i)] We have $|A \cdot A \cdot A| \lesssim K^{O(1)} |A|$.
\item[(ii)] There exists a $O(K^{O(1)})$-approximate group $\tilde A$ of size $|\tilde A| \sim K^{O(1)} |A|$ 
which contains $A$.
\item[(iii)] There exist $O(K^{O(1)})$-approximate groups $B_1 \subseteq B_2 \subseteq B_3 \subset H$ and 
$C \subset G/H$ with 
\begin{equation}\label{b-compare}
|B_3| \lesssim K^{O(1)} |B_1|, 
\end{equation} 
together with a partial inverse $\phi: C^3 \to G$ to $\pi$, with $\phi(1) = 1$ and $\phi(x^{-1}) = \phi(x)^{-1}$ for all $x \in C^3$,
such that
\begin{equation}\label{phib-approx}
 \phi(x) B_i  \subseteq B_{i+1} \phi(x); \quad
 B_i \phi(x)  \subseteq \phi(x) B_{i+1} \quad
  \hbox{ for all } x \in C, i = 1,2
\end{equation}
and
\begin{equation}\label{phib-2-approx}
\phi(x) \phi(y) \phi(z) \in \phi(xyz) B_3 \hbox{ for all } x,y,z \in C
\end{equation}
with the containment
\begin{equation}\label{adef-approx}
 A \subseteq \bigcup_{x \in C} \phi(x) B_1
\end{equation}
and the cardinality bound
\begin{equation}\label{acard-lower}
|A| \gtrsim K^{-O(1)} |B_1| |C|.
\end{equation}
\end{itemize}
\end{lemma}

\begin{remark} One can extend this lemma to the continuous setting provided that the topology of $H$ is well-behaved (e.g. the projection map $\pi$ should be continuous and open, and in particular $H$ should be closed) and one can ``disintegrate'' the measure $\mu$ on $G$ into the measures on $H$-cosets of $G$, integrated against the measures on $G/H$; this can for instance be done if $G$ is a finite-dimensional Lie group 
and $H$ is a closed Lie subgroup. We omit the details.  There is likely to also be an entropy analogue of this lemma under reasonable assumptions on the metric but we will not describe these here.  The approximate homomorphism $\phi$ is only defined on $C^3$, but an inspection of the proof below shows that one could in fact extend it to $C^n$ for any fixed $n$, and have a sequence $B_1 \subseteq \ldots \subseteq B_n$ of nested approximate groups of comparable size, with suitable modifications to \eqref{b-compare}, \eqref{phib-approx}, \eqref{phib-2-approx}; again, we omit the details.
\end{remark}

\begin{proof} The equivalence of (i) and (ii) follows from Theorem \ref{tripling-classify}.  Now let us see that (iii) implies (i).  From \eqref{adef-approx} we have
$$ A \cdot A \cdot A \subseteq \bigcup_{x,y,z \in C} \phi(x) B_1 \phi(x_2) B_1 \phi(x_3) B_1.$$
From repeated application of \eqref{phib-approx} we have
$$ \phi(x) B_1 \phi(y) B_1 \phi(z) B_1 \subseteq \phi(x) \phi(y) \phi(z) B_3 B_2 B_1$$
and hence by \eqref{phib-2-approx}
$$ A \cdot A \cdot A \subseteq \bigcup_{x,y,z \in C} \phi(xyz) B_3 B_3 B_2 B_1 \subseteq \bigcup_{w \in C^3} \phi(w) B_3^4$$
and thus $|A^3| \leq |C^3| |B_3^4|$.  The claim now follows from \eqref{acard-lower}, \eqref{b-compare}, and the hypothesis that $C$ and $B_3$
are $O(K^{O(1)})$-approximate groups.

It remains to show that (ii) implies (iii).  By replacing $A$ by $\tilde A$ if necessary we may assume that $A$ is itself a $O(K^{O(1)})$-approximate group; in particular, $A$ is symmetric, contains $1$, and (by Theorem \ref{tripling-classify})
we have $|A^n| \lesssim_n K^{O_n(1)} |A|$ for all $n \geq 1$.  
Since $\pi$ is a homomorphism, we easily see that $\pi(A)$ is also a $O(K^{O(1)})$-approximate group.  Thus we shall set $C := \pi(A)$.
Now we construct the $B_i$ by the formulae
$$ B_1 := (A^2 \cap H)^3; \quad B_2 := (A^8 \cap H)^3; \quad B_3 := (A^{26} \cap H)^3.$$
Observe that if $a, a' \in A$ lie in the same fiber of $\pi$ (i.e. in the same coset of $H$), then $a' \in (A^2 \cap H) a$.  Since $A$ intersects exactly $|C|$ fibers of $\pi$, we conclude that $|A| \leq |C| |A^2 \cap H|$.  On the other hand, observe that if $a \in A$, then the set $(A^{2n} \cap H) a$ lies in $A^{2n+1}$, and is also contained in the same fiber of $\pi$ as $a$.  This implies that $|A^{2n+1}| \geq |C| |A^{2n} \cap H|$ for all $n \geq 1$.  Since $|A^{2n+1}| \lesssim_n K^{O_n(1)} |A|$, we conclude
that $|A^{2n} \cap H| \sim_n K^{O_n(1)} |A^2 \cap H|$ for all $n \geq 1$.  Since $(A^{2n} \cap H)^3 \subseteq A^{6n} \cap H$, we conclude
that $A^{2n} \cap H$ has a tripling constant of $O_n(K^{O_n(1)})$.
Applying Corollary \ref{tripling-better} (observing that $A^{2n} \cap H$ is symmetric and contains $1$)
we thus see that $B_1, B_2, B_3$ are all $O(K^{O(1)})$-approximate groups. Note that the estimates established here also give \eqref{acard-lower} and
\eqref{b-compare}.

Since $C = \pi(A)$, we have $C^3 = \pi(A^3)$.  Using the axiom of choice (which is actually unnecessary here, since $A^3$ and $C^3$ are finite sets), we may select a partial inverse $\phi: C^3 \to A^3$ to $\pi$, which takes values in $C$ on $A$; since $C$ and $A$ are both symmetric and contain the origin, there is no difficulty requiring $\phi(1) = 1$ and $\phi(x^{-1}) = \phi(x)^{-1}$ for all $x \in C^3$.  If $a$ and $\phi(x)$ lie in the same fiber of $\pi$, then as mentioned before we have $a \in \phi(x) (A^2 \cap H)$, which implies \eqref{adef-approx}.

Now observe that if $x \in C$, then $\phi(x) (A^{2n} \cap H)^3 \phi(x)^{-1}$ lies in $A^{6n+2} \cap H$, and hence
$\phi(x) B_i \phi(x)^{-1}$ lies in $B_{i+1}$ for $i=1,2$.  This proves the inclusions $\phi(x) B_i \subseteq \phi(x) B_{i+1}$; the reverse inclusions in
\eqref{phib-approx} are proven similarly.

Finally, for $x,y,z \in C$, we observe that $\phi(x) \phi(y) \phi(z) \phi(xyz)^{-1} \in A^4 \cap H$, and so \eqref{phib-2-approx} follows.
This concludes the implication of (iii) from (ii).
\end{proof}

In principle, this splitting lemma should allow one to deduce inverse product set estimates for nilpotent (and perhaps even solvable) groups
from the abelian theory.  To do this in full generality appears to be rather difficult however, and we shall only demonstrate the situation with
a particularly simple nilpotent group, namely a Heisenberg group.

\begin{definition}[Heisenberg group]  Let $Z, W$ be additive abelian groups, and let $\{,\}: Z \times Z \to W$ be an antisymmetric mapping (thus $\{x,y\} = -\{y,x\}$) which is a homomorphism in each of the two variables separately (thus $\{ x+y, z \} = \{x,z\} + \{y,z\}$ and $\{x,y+z\} = \{x,y\} + \{x,z\}$).  We define the \emph{Heisenberg group} associated to this antisymmetric mapping to be the set $G := Z \times W$ endowed with the group law
$$ (z, w) \cdot (z',w') := (z + z', w + w' + \{ z, z' \} ).$$
\end{definition}

One easily verifies that $G$ is a discrete multiplicative group with the ``vertical'' group $\{0\} \times W$ (which we identify with $W$) as a normal subgroup (indeed, it lies in the centre of $G$) and with identity $1_G = (0,0)$ and inverse $(z,w)^{-1} = (-z,-w)$;
the quotient $G/W$ is canonically identified with the ``horizontal'' group $Z$. Thus we have the short exact sequence
$$ 0 \to W \to G \to Z \to 0.$$
Since $W$ and $Z$ are abelian, we thus see that $G$ is a $2$-step nilpotent group.

We write $Z \times W$ for the additive group which is the product of $Z$ and $W$, and $\iota: G \to Z \times W$ for the identity map from $G$
to $Z \times W$.  We caution that while $\iota$ is a bijection, it is \emph{not} a group homomorphism from the multiplicative group $G$ to the additive group $Z \times W$.  Nevertheless, $\iota$ does relate the subgroups of $G$ with the subgroups of $Z \times W$ (except for a ``$2$-torsion'' issue)
as follows.  Let $\pi: Z \times W \to Z$ be the canonical projection, and for any $A, B \subset Z$ let $\{A,B\} \subset W \equiv \{0\} \times W$
denote the set $\{A,B\} := \{ \{a,b\}: a \in A, b \in B \}$.  We also let $\langle \{A,B\} \rangle$ denote the subgroup of $\{0\} \times W$ 
generated by $\{A,B\}$.  Finally, given any $A \subset Z \times W$ we write $2 \cdot A := \{ 2x: x \in A \} = \{ (z+z, w+w): (z,w) \in A \}$.

\begin{proposition}[Subgroups of the Heisenberg group]  
Let $G$ be a Heisenberg group arising from an antisymmetric mapping $\{,\}: Z \times Z \to W$, and let 
$A \subset G$ be a multiplicative subgroup of $G$.  Then there exists an additive subgroup $\tilde A$ of $Z \times W$ such that
\begin{equation}\label{binclude}
 2 \cdot(\tilde A + \langle \{ \pi(\tilde A), \pi(\tilde A) \} \rangle ) \subseteq \iota(A) \subseteq \tilde A + \langle \{ \pi(\tilde A), \pi(\tilde A) \} \rangle 
 \end{equation}
\end{proposition}

\begin{remark} In the converse direction, it is easy to verify that given any additive subgroup $\tilde A$ of $Z \times W$, that the set $\iota^{-1}( \tilde A + \langle \{ \pi(\tilde A), \pi(\tilde A) \} \rangle )$ is a multiplicative subgroup of $G$.  Thus the above proposition classifies multiplicative subgroups of $G$ in terms of additive subgroups of
$Z \times W$, except for the ``$2$-torsion'' issue of having to distinguish a set the additive group $A' := \tilde A + \langle \{ \pi(\tilde A), \pi(\tilde A) \} \rangle$ from its dilate $2 \cdot A'$.  If $Z \times W$ is finitely generated, then the quotient of the additive group $A'$ by $2 \cdot A'$ will be bounded, and so in some sense the above classification of subgroups of $G$ only ``loses'' a bounded amount of information.
\end{remark}

\begin{proof}  First observe that if $(z,w)$ and $(z',w')$ are in $A$, then 
$$ (z,w) \cdot (z',w') \cdot (z,w)^{-1} \cdot (z',w')^{-1} = (0, 2 \{z,z'\} ).$$
Thus if we set $C := \pi(\iota(A)) \subset Z$, we see that $2 \cdot \{C,C\} \subset A$, and hence (since the Heisenberg group law is additive on
$\{0\} \times W$) we also have $2 \cdot \langle \{C,C\} \rangle \subset A$.  If we now set $\tilde A := \iota(A) + \langle \{C,C\} \rangle$, we easily verify that $\tilde A$ is indeed an additive group and obeys the inclusions \eqref{binclude}.
\end{proof}

We now extend this proposition to approximate groups, though to deal with the $2$-torsion issue we shall need to make an additional assumption on the ``vertical'' group $W$.

\begin{theorem}[Approximate subgroups of the Heisenberg group]\label{heisen}  
Let $G$ be a Heisenberg group arising from an antisymmetric mapping $\{,\}: Z \times Z \to W$ such that $W$ has no $2$-torsion (thus if $w \in W$ and $2 w = 0$, then $w=0$), and let $A \subset G$ be a finite nonempty subset of $G$ such that $|A \cdot A \cdot A| \leq K |A|$ for some $K \geq 1$.
Then there exists a $O(K^{O(1)})$-approximate additive subgroup $\tilde A$ of $Z \times W$, such that
\begin{equation}\label{pi-include}
\{ \pi( \tilde A ), \pi(\tilde A) \} \subseteq \tilde A
\end{equation}
and
$$ \iota(A) \subseteq \tilde A.$$
Furthermore we have
$$ |A| \gtrsim K^{-O(1)} |\tilde A|.$$
\end{theorem}

\begin{remark} In the converse direction, if $\tilde A$ obeys all the above properties, then a comparison of the multiplicative group law for $G$
and the additive group law for $Z \times W$ reveals that
$$ \iota(A \cdot A \cdot A) \subseteq 3 \tilde A + 3 \{ \pi(\tilde A), \pi(\tilde A) \} \subseteq 6 \tilde A$$
and hence by the approximate group properties of $\tilde A$ 
$$ |A \cdot A \cdot A| \leq 6 |\tilde A| \lesssim K^{O(1)} |\tilde A| \lesssim K^{O(1)} |A|.$$
Thus we have a sharp characterisation of the sets of small tripling in the Heisenberg group $G$, in the case when no $2$-torsion is present in the centre.  In principle, one can make this characterisation more explicit by using a version of Freiman's theorem (such as the one in \cite{gr-4})
in the abelian group $Z \times W$ to classify $\tilde A$, and then to work with the concrete description of $\tilde A$ given by that 
theorem to determine which approximate groups $\tilde A$ obey the constraint \eqref{pi-include}.  Of course once one characterises sets of small tripling, one can use the results of earlier sections to characterise sets of small doubling, or of with small partial product set, etc.  These fully explicit descriptions are however rather lengthy to state and we will leave them to the reader.
\end{remark}

\begin{proof}  We apply Lemma \ref{split-approx} with $H := \{0\} \times W \equiv W$ to obtain
$O(K^{O(1)})$-approximate groups\footnote{It is a somewhat unfortunate circumstance that we will be regarding $Z$ and $W$ both as additive groups and multiplicative group (with the same group operation); thus for instance if $B \subset W$ then $B+B = B \cdot B$.  We hope the reader will not be unduly confused by this.}
 $B_1 \subseteq B_2 \subseteq B_3 \subset W$ and $C \subset Z$, and a partial inverse $\phi: C^3 \to G$ to the projection map $\pi: G \to Z$
with $\phi(0) = (0,0)$ and $\phi(-x) = -\phi(x)$ that obeys \eqref{b-compare}, \eqref{phib-2-approx}, \eqref{adef-approx}, \eqref{acard-lower}.  (Because $W$ is abelian, the containments \eqref{phib-approx} become trivial and will not be needed here.)   Since $\phi$ is a partial inverse to $\pi$, we may 
 write $\phi(z) = (z, f(z))$ for some odd function $f: C^3 \to W$.  Thus for instance \eqref{adef-approx} becomes the assertion that
 $w \in f(z) + B_1$ for all $(z,w) \in A$.  From \eqref{phib-2-approx} (setting the third element of $C$ to be the identity) we see that
 for all $z_1,z_2 \in C$ we have
$$ (z_1, f(z_1)) \cdot (z_2, f(z_2)) \in (z_1+z_2, f(z_1+z_2)) \cdot B_3$$
and hence (expanding out the group multiplication law in coordinates)
$$ f(z_1) + f(z_2) + \{ z_1, z_2 \} \in f(z_1+z_2) + B_3.$$
Swapping $z_1$ and $z_2$ and then subtracting, we see that
$$ 2 \cdot \{ z_1, z_2 \} \in B_3 - B_3 \hbox{ for all } z_1, z_2 \in C.$$
Let $\tilde B := (B_3 - B_3) \cap (2 \cdot W)$.  The set $2 \cdot W$ is a subgroup of the abelian group $W$, and so
$\tilde B + \tilde B + \tilde B = (3B_3 - 3B_3) \cap (2 \cdot W)$.  Since $B_3$ is a $O(K^{O(1)}$-approximate group, we can cover $3B_3 - 3B_3$
by $O(K^{O(1)})$ translates of $B_3$, and hence $\tilde B + \tilde B + \tilde B$ can be covered by $O(K^{O(1)})$ translates of $B_3$, intersected
with $2 \cdot W$.  Any one of these translates is itself contained in a translate of $(B_3-B_3) \cap (2 \cdot W)$, and so we see that
$|3 \tilde B| \lesssim K^{O(1)} |\tilde B|$.  Applying Lemma \ref{tripling-better} we see that $3\tilde B$ is a $O(K^{O(1)})$-approximate group.
If we then let
$$ B' := \{ b \in W: 2 \cdot b \in 3 \tilde B \}$$
we conclude (since the map $b \mapsto 2 \cdot b$ is a group isomorphism between $W$ and $2 \cdot W$) that $B'$ is also a $O(K^{O(1)})$-approximate group.  By construction we have
\begin{equation}\label{fuzz}
 \{ C, C \} \subseteq B'.
\end{equation}
Next, observe that
$$ |B' + 3\tilde B| = |B' + 2 \cdot B'| \leq |3 B'| \lesssim K^{O(1)} |B'| = K^{O(1)} |3\tilde B|$$
and hence by Lemma \ref{cover}, we can cover $B'$ by $O(K^{O(1)})$ translates of $3\tilde B - 3\tilde B$, which is contained in $12 B_3$.  Since $B_3$ is itself a $O(K^{O(1)})$-approximate group, we can conclude that
\begin{equation}\label{namby}
|nB' + mB_3| \lesssim_{n,m} K^{O_{n,m}(1)} |B_3|
\end{equation}
for any $n,m \geq 1$.  To use this, define the set $A' \subset Z \times W$ by
$$ A' := \{ (z, w): z \in C, w \in f(z) + 9B' + B_1 \}.$$
From \eqref{adef-approx} we see that $A'$ contains $A$.  Also we have
\begin{align*}
3\iota(A') &= \bigcup_{z_1,z_2,z_3 \in C} (z_1,f(z_1) + 9B' + B_1) + (z_2,f(z_2) + 9B' + B_1) + (z_3,f(z_3) + 9B' + B_1) \\
&\subseteq \bigcup_{z_1,z_2,z_3 \in C} (z_1+z_2+z_3, f(z_1)+f(z_2)+f(z_3)+ 27B' + 3B_3).
\end{align*}
On the other hand, from \eqref{phib-2-approx} we have
$$ f(z_1)+f(z_2)+f(z_3) + \{ z_1, z_2 \} + \{ z_1, z_3 \} + \{ z_2, z_3 \} \in f(z_1+z_2+z_3) + B_3$$
and hence by \eqref{fuzz}
$$ f(z_1)+f(z_2)+f(z_3) \in f(z_1+z_2+z_3) + B_3 + 3B'.$$
Thus we have
$$ 3 \iota(A') \subseteq \bigcup_{z_1,z_2,z_3 \in C} (z_1+z_2+z_3, f(z_1+z_2+z_3)+ 30B' + 4B_3)$$
and hence
$$ |3 \iota(A')| \leq |C^3| |30B'+4B_3|.$$
Since $C$ is a $O(K^{O(1)})$-approximate group, we thus see from \eqref{namby} and \eqref{acard-lower} that
$$ |3\iota(A')| \lesssim K^{O(1)} |A|.$$
Since $|\iota(A')| \geq |A|$, we thus see from Lemma \ref{tripling-better} that the set $\tilde A := 3 (\iota(A') \cup 0 \cup - \iota(A'))$ is
a $K^{O(1)}$-approximate group with
$$ |\tilde A| \sim K^{O(1)} |A|.$$ 
Finally, we see from construction that $\pi(\tilde A) \subseteq 3C$, and hence $\{ \pi(\tilde A), \pi(\tilde A) \} \subseteq 9 \{C,C\} \subseteq
9B'$ by \eqref{fuzz}; since $f(0) = 0$, we obtain \eqref{pi-include} as desired.
\end{proof}

\end{document}